\documentclass[letter,9pt]{article}
\usepackage[scale=.83]{geometry}
\usepackage[affil-it]{authblk}
\usepackage{amssymb,amsmath,amsthm,amsbsy}

\usepackage{physics}
\usepackage{graphicx}
\usepackage{epstopdf}
\usepackage{epsfig}
\usepackage{multicol}
\usepackage{graphicx}
\usepackage{caption}
\usepackage{subcaption}
\usepackage{multirow}
\usepackage{algorithm,algorithmic}
\usepackage{setspace}
\usepackage{color}
\usepackage{tablefootnote}
\usepackage{nomencl}
\makenomenclature

\usepackage{lineno}
\usepackage{svg}
\usepackage{pdfpages}
\usepackage{comment}
\usepackage{booktabs}
\usepackage{nomencl}
\makenomenclature
\usepackage{ifthen}

\usepackage{microtype} % Slightly tweak font spacing for aesthetics

\usepackage{titlesec} % Allows customization of titles
\usepackage[toc]{appendix}
\usepackage{hyperref}
\usepackage{cleveref}

\titleformat{\section}[block]{\large\scshape\centering}{\thesection.}{1em}{} % Change the look of the section titles
\titleformat{\subsection}[block]{\large}{\thesubsection.}{1em}{} % Change the look of the section titles

\title{\textbf{Data-Driven Reduction of the Finite-Element Model of a Tribomechadynamics Benchmark Problem}}

\author{Ahmed Amr Morsy\footnote{Corresponding author: morsya@ethz.ch}\ } \author{Zhenwei Xu}
\author{Paolo Tiso}
\author{George Haller} 
\affil{
	Institute for Mechanical Systems, ETH Zürich, Leonhardstrasse 21, 8092 Zurich, Switzerland}

\date{}

\renewcommand{\vec}[1]{\boldsymbol{\mathrm{#1}}}     
\newcommand{\mat}[1]{\boldsymbol{#1}}           % for matrices

\newcommand{\mt}[3]{\mat{#1}_\text{#2}^\text{#3}}
\newcommand{\vt}[3]{\vec{#1}_\text{#2}^\text{#3}}
\newcommand{\st}[3]{\text{#1}_\text{#2}^\text{#3}}

\begin{document}
	
\graphicspath{{figures/}}

\maketitle 
\thispagestyle{empty}

    	\section{abstract} 
     Bolted joints can exhibit nonsmooth and significantly nonlinear dynamics. Finite Element Models (FEMs) of this phenomenon require fine spatial discretizations, inclusion of nonlinear contact and friction laws, as well as geometric nonlinearity. Owing to the nonlinearity and high dimensionality of such models, full-order dynamic simulations are computationally expensive. In this work, we use the theory of Spectral Submanifolds (SSMs) to construct a data-driven, smoothed reduced model for a 187,920-dimensional FEM model of a broadly studied  Tribomechadynamics benchmark structure with bolted joints.  We train the 4-dimensional reduced model using only a few transient trajectories of the full unforced FEM model. We show that this smooth model accurately predicts the experimentally observed nonlinear forced response of the full nonsmooth benchmark problem.

        \noindent\textbf{Keywords}: joints, spectral submanifolds, nonlinear dynamics, model order reduction, geometric nonlinearity\\	

    \section{Introduction}
\label{section:introduction}
Microslip refers to the partial slipping between the surfaces of the components brought in contact. This is in contrast to complete slip, where all the resistance due to friction is lost, and in contrast to complete stick, where there structure is linear. Microslip results in a degrading amplitude-dependent stiffness and, more importantly, a nonlinear damping that becomes the main source of dissipation \cite{Brake2018b}. Predictive numerical models of this phenomenon typically require High-Fidelity (HF) Finite Element Models (FEMs), where a fine mesh of frictional contact elements on the interfaces is used. Time-integration of these nonlinear, high-dimensional equations of motion is computationally expensive. This has generated a substantial effort to derive Reduced Order Models (ROMs) for joints.

\begin{figure}
    \centering
    \includegraphics[scale=0.32]{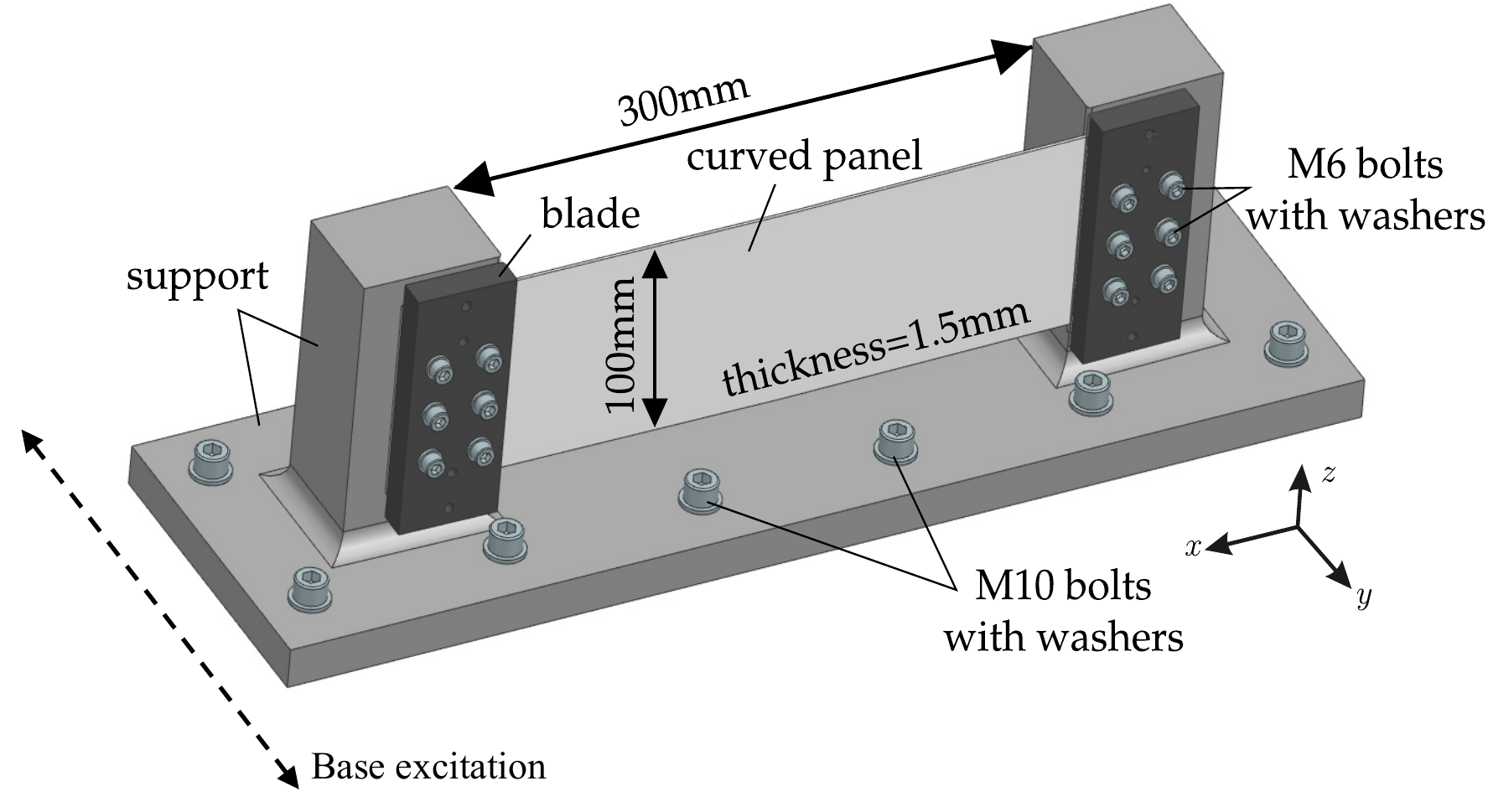}
    \caption{CAD Model of TRChallenge benchmark structure.\cite{darus-3147_2022}}
    \label{fig:TRC_CAD}
\end{figure}

As a survey on the state-of-the-art ROM methods, we recall the recent 2021 Tribomechadynamics Research Challenge (TRChallenge) \cite{TRCarticle}, in which eight international research groups participated. The TRChallenge involved making blind predictions of the forced nonlinear dynamics of a mechanical assembly that is yet to be fabricated and tested. The assembly, shown in Fig. \ref{fig:TRC_CAD}, consists of a slender panel that is preclamped to the support structure using bolts. This design triggers a geometric nonlinear response due to the bending-stretching coupling of the thin panel, in addition to friction nonlinearity due to the slippage of the panel at the joints. The forcing is a base excitation that triggers the first nonlinear modal response featuring the out-of-plane bending of the curved panel. The required modal characteristics were the amplitude-dependent resonance frequencies and modal damping. Experimentally, the response additionally featured an internal resonance with the torsional mode of the panel. 

Krack et al. \cite{TRCarticle} documented the details and the results of the eight prediction approaches used. Only five of these approaches included geometric nonlinearity in the model. Among those, one approach relied on post-processing the signal of a full transient simulation of a FEM. Full simulations, however, cannot be used to predict the forced dynamics of FEMs due to the associated prohibitive costs. Moreover, in the case of internally resonant dynamics, as in the case of this structure, extracting modal characteristics from a decaying signal is not straightforward. The remaining four approaches constructed ROMs. One approach used Quasi-Static Modal Analysis (QSMA) \cite{Lacayo2019b}. QSMA solves a set of nonlinear static problems due to an external static load having the mode shape of the structure. As an analysis method, it is derived from the modal analysis of linear, undamped systems. Therefore, while it is true that QSMA may, in some cases, provide useful preliminary insight into the softening/hardening behavior of weakly nonlinear joints, it cannot fully replace dynamic analyses. Multi-harmonic responses, velocity-dependent forces, nonlinear dynamics not restricted to low-dimensional linear subspaces, internal resonances, and external forcing are all circumstances under which QSMA will not be useful.

Another ROM featured in \cite{TRCarticle} assumed a single-degree-of-freedom model with mono-harmonic response, and assumes a decoupling of geometric and friction nonlinearities which are then analyzed separately. Despite computing good predictions of the frequencies, this method cannot be generalized due to the aforementioned assumptions. The remaining two approaches in \cite{TRCarticle} implemented a substructuring step \cite{Allen} to reduce their respective FEMs. Substructuring aims to represent the displacement field based on the mode shapes of individual components in addition to static deformations along their shared boundaries. There is no guarantee, however, that this is a valid representation. For instance, in one of the two approaches that used substructuring, the overly stiff response was attributed to retaining an insufficient number of component modes. As for the second approach, its overly stiff results partly arise from assuming constant normal stresses on the interface of the joint. Notably, a model-driven ROM technique for jointed structures, based on hyperreduction, was proposed in \cite{Morsy2023}. While efficient, the accuracy of the resulting ROM depends on the projection basis, whose construction can be challenging. In light of the difficulties discussed, we propose creating HF FEMs with no assumptions on the dynamic response or on the behavior of the joint. Then, instead of developing an equation-driven ROM, we use a data-driven approach based on the theory of Spectral Submanifolds (SSMs).

The recent theory of SSMs proves the existence and persistence of attracting, invariant manifolds of smooth nonlinear systems under nonresonance conditions \cite{Haller2016}. This result has been used for constructing reduced-order models (ROMs) of mechanical systems. The numerical open-source package SSMTool has been developed for model order reduction of high-dimensional mechanical systems \cite{Ponsioen2020,Jain2022}, that can also exhibit internal resonances \cite{Li2022,Li2022b}, and parametric resonances \cite{Thurnher2024}. In a data-driven setting, SSMLearn \cite{Cenedese2022a,Cenedese2022,SSMLearn} and FastSSM \cite{Axas2022} have been developed based on the theory of SSMs. Data-driven ROMS have been created to successfully predict dynamic responses in fluid problems \cite{Kaszas2022, Haller2023}, control problems \cite{Alora2023}, and structural nonlinear dynamics \cite{Cenedese2022a,Cenedese2022,Axas2022}.

Recently, Cenedese et al. \cite{Cenedese2023} used SSMLearn for data-driven reduced order modeling of high-dimensional smooth FE models such as internally resonating geometrically nonlinear structures. In this work, we adopt this methodology and apply it to a FEM of the 2021 Tribomechadynamics benchmark discussed extensively by Krack et al. \cite{TRCarticle}. Even though the full numerical model is nonsmooth, its flow map will be at least continuous, hence its invariant manifolds are also expected to be at least continuous. By performing data-driven SSM reduction to trajectories of this full model, we effectively approximate the dominant continuous SSM with the closest fitting smooth SSM. This approach is non-intrusive since the SSM-reduced nonlinear dynamics are learned directly from data. There exist other numerical toolboxes, such as SINDy \cite{SINDy}, that fit sparse models to low-dimensional data. Our SSM-based ROM, however, is guaranteed to persist under time-dependent forcing \cite{Haller2016,Haller2024,xu2024}. This enables us to test the predictive power of our approach by how well the smooth SSM-based model (trained only on unforced data) is able to predict the nonsmooth forced response. As we will see from a detailed comparison with \cite{TRCarticle}, our SSM-based prediction is highly accurate and surpasses all available results from other methods.

In the next section, we review some results of the theory Spectral Submanifolds (SSMs) for smooth mechanical systems. In section \ref{section:joints}, we set up the FEM of the joint and explain our data-driven reduced order modeling approach using the numerical package SSMLearn. Next, we illustrate our approach in Section \ref{section:TRCBenchmark} by showing results for a high-fidelity FEM of the TRCBenchmark structure. Lastly, we present our conclusions in Section \ref{section:conclusions}.

    \section{Background}    \label{section:Background}
         Consider the equations of motion of a FEM of a smooth mechanical system
    \begin{equation}
    \mat{M}\ddot{\vec{q}} + \mat{C}\dot{\vec{q}} + \mat{K}\vec{q} + \vt{f}{}{int} (\vec{q},\dot{\vec{q}}) = \epsilon \vt{f}{}{ext}(\vec{\Omega} t),
    \label{eq:EqnsMotion2ndOrder}
    \end{equation}
    \noindent
    where $\vec{q}$, $\dot{\vec{q}}$, $\ddot{\vec{q}} \in \mathbb{R}^{n}$ are the generalized nodal displacement, velocity, and acceleration vectors with $n$ being the number of degrees of freedom of the system. The matrices $\mat{M}$, $\mat{C}$, $\mat{K} \in \mathbb{R}^{n \cross n}$ are the mass, damping and linear stiffness matrices of the unconnected structure, $\vt{f}{}{int}$ is a vector of nonlinear forces, which here features contact, friction, and geometrically nonlinear forces, $\vt{f}{}{ext} \in \mathbb{R}^{n}$ is the vector of external forces imposed on the structure with frequencies of excitation contained in the vector $\vec{\Omega} \in \mathbb{R}^{k}$ and a small parameter $\epsilon>0$ scaling the forcing amplitude. 
    \noindent
    We introduce $\vec{x}=[\vec{q},  \vec{\dot{q}}]^{\text{T}} \in \mathbb{R}^{2n}$ to express Eq.\eqref{eq:EqnsMotion2ndOrder} in $1^{\text{st}}$ order form as

    \begin{equation}
        \vec{\dot{x}} =
        \mat{A} \vec{x} + \vt{f}{0}{}(\vec{x}) + \epsilon \vt{f}{1}{}(\vec{\Omega}t),
        \label{eq:smooth_1st_order}
        \\
        \end{equation}
        \begin{equation*}
        \mat{A} =
        \begin{bmatrix}
            \mat{0} && \mat{I}\\
            -\mt{M}{}{-1}\mt{K}{}{} &&  -\mt{M}{}{-1}\mat{C}
        \end{bmatrix}, \qquad
        \vt{f}{0}{} =
        \begin{bmatrix}
            \vec{0}\\
            -\mt{M}{}{-1}\vt{f}{}{int}
        \end{bmatrix}, \qquad
        \vt{f}{1}{} =
        \begin{bmatrix}
            \vec{0}\\
            \mt{M}{}{-1}\vt{f}{}{ext}(\vec{\Omega} t)
        \end{bmatrix}.               
        \end{equation*}
    \noindent
    We assume that $\vec{x}=\vec{0}$ is a fixed point of the $\epsilon=0$ limit of system \eqref{eq:smooth_1st_order}.

    \subsection{Linear Dynamics}
    Consider the linear dynamics
    \begin{equation}
        \vec{\dot{x}} = \mat{A} \vec{x}
        \label{eq:lineardynamics}
    \end{equation}
    in Eq. \eqref{eq:smooth_1st_order}. Since we are studying oscillatory structural dynamics, we focus on the case where the modes of the system are underdamped. We write the eigenvalue problem associated with $\mat{A}$ as $\mat{A}\mat{\Phi}{}{} = \mat{\Lambda} \mat{\Phi}{}{}$, where $\mat{\Lambda}$ is a diagonal matrix with complex conjugate eigenvalues $\lambda_{1},\bar{\lambda}_{1}, \lambda_{2},\bar{\lambda}_{2}, ...,\bar{\lambda}_{n}$ on the diagonal, ordered such that
    \begin{equation}
        \text{Re}\ \lambda_{n} \leq \text{Re}\ \lambda_{n-1} \leq ... \leq \text{Re}\ \lambda_{2} \leq \text{Re}\ \lambda_{1} < 0.
        \label{eq:realparts}
    \end{equation}
    The columns of $\mat{\Phi} \in \mathbb{C}^{2n\cross2n}$ are the corresponding eigenvectors ${\vec{\phi}_{1}, \bar{\vec{\phi}}_{1}, \vec{\phi}_{2}, \bar{\vec{\phi}}_{2}, ..., \bar{\vec{\phi}}_{n}} \in \mathbb{C}^{2n}$. We express the eigenvalues $\lambda_{j}, \bar{\lambda}_{j}$ of the $j$-th mode of the system as $\alpha_{j} \pm i \omega_{j}$. We have assumed in \eqref{eq:realparts} that all eigenvalues have negative real parts, i.e. the linear system is asymptotically stable, and that the columns of $\mat{\Phi}$ are linearly independent. The real part $\alpha_{j}$ describes the decay exponent of the $j$-th mode, while $\omega_{j}$ is its frequency. The corresponding eigenvectors  $\vec{\phi}_{j}, \bar{\vec{\phi}}_{j}$ form a 2-dimensional linear subspace. For a real representation of this spectral subspace, we define $E_{j} = \{\text{Re}\{\vec{\phi}_{j}\} \cross \text{Im}\{\vec{\phi}_{j}\}\} \subset \mathbb{R}^{2n}$. The spectral subspace $E_{j}$ is an invariant subspace of the linear system, meaning that a trajectory $\vec{x}(t)$ having initial conditions $\vec{x}(t_{0}) \in E_{j}$ remains in $E_{j}$. The slowest $m$ modes of system \eqref{eq:lineardynamics} are those with the decay exponents  $\alpha_{j}$. Accordingly, we define $E^{m}$ as the linear subspace formed by the slowest $m$ modes,  $E^{m} := \{\st{E}{1}{} \cross \st{E}{2}{} \cross ... \cross E_{m}\}$. This linear subspace is attracting, meaning that trajectories starting from generic initial conditions in the neighborhood of the fixed point eventually get arbitrarily close to $E^{m}$. This attracting property of $E^{m}$, together with its invariance, make $E^{m}$ an ideal candidate for model reduction in the linear systems \eqref{eq:lineardynamics}.

    \subsection{Nonlinear Dynamics}
    \label{subsection:NonlinearDynamics}

    Consider now the nonlinear, smooth autonomous system 

    \begin{equation}
        \vec{\dot{x}} =\mat{A} \vec{x} +\vt{f}{0}{}(\vec{x}).
        \label{eq:NLdynamics}
    \end{equation}
    Under a non-resonance assumption\footnote{Namely, the spectral subspace $E^{2m}$ of the semi-simple matrix $\mat{A}$ is non-resonant. (i.e., no nonnegative, low-order, integer linear combination of the spectrum inside $E^{2m}$ is contained in the spectrum of $\mat{A}$ outside $E^{m}$) \cite{Haller2016}.This condition can always be fulfilled by increasing the dimensionality of the spectral subspaces $E^{2m}$, and accordingly, the dimension of the spectral submanifolds.} on the eigenvalues of $\mat{A}$, by the theory of Spectral Submanifolds (SSMs) \cite{Haller2016}, dynamics do not generally occur in a linear spectral subspace $E^{m}$ anymore, but in perturbations of $E^{m}$ that are nonlinear continuations of it. These nonlinear continuations of the spectral subspaces are called spectral submanifolds, and are all tangent to $E^{m}$ at the fixed point $\vec{x}=\vec{0}$. Each of these submanifolds contains local solutions of the nonlinear system, making them the geometrical entities of interest since invariance is fulfilled on them. Further, these submanifolds have the same $m$-dimensionality as $E^{m}$, inherit the property of being attracting, and continue to exist even if small perturbations are introduced to the vector field \cite{Haller2016}. Interestingly, among all those submanifolds, a single smoothest one, referred to as the (primary) SSM \cite{Haller2016,Haller2023} exists. We denote the primary SSM by $\mathcal{W}(E^{m})$. All the other spectral submanifolds that are less smooth than the primary SSM are called secondary SSMs \cite{Haller2023}. The smoothness property of $\mathcal{W}(E^{m})$ is leveraged for computations. In particular, one can construct an invariant, robust Reduced Order Model (ROM) using a Taylor expansion (\cite{Haller2016,Ponsioen2020,Jain2022}). For $\epsilon > 0$, and under time-dependent forcing with $k>0$, under one additional assumption\footnote{Namely, no nonnegative, low-order, integer linear combination of the real parts of the eigenvalues of $\mat{A}$ inside $E^{m}$ coincides with the real part of the spectrum of any eigenvalue of $\mat{A}$ outside $E^{m}$ \cite{Haller2016}.}, the submanifolds become time-dependent, and the time-dependent primary SSM remains ideal for model reduction.

        \section{Joints}
    \label{section:joints}
    In this section, we consider FE models of jointed mechanical assemblies. We focus on structures exhibiting microslip behavior \cite{Gaul1997}. This refers to a kinematic regime in which partial slipping occurs between interfaces in contact. Loads are thus transmitted between the joined components via friction while exhibiting pronounced nonlinear damping. Microslip can be modeled by employing a fine mesh of frictional contact elements on the interfaces \cite{Krack2017}. In this work, we use 3D frictional contact elements, readily available in many commercial finite element software. Tangentially, friction is defined using a planar Coulomb law. If the nodal pairs are in contact, the tangential slip forces are 
    
  \begin{gather}
        f_{t,u} =
          \mu f_{n} \frac{\dot{u}}{\sqrt{\dot{u}^{2}+\dot{v}^{2}}}, 
         \qquad
        f_{t,v} =
          \mu f_{n} \frac{\dot{v}}{\sqrt{\dot{u}^{2}+\dot{v}^{2}}} \qquad \text{if } 
        \dot{u} \neq 0\ ||\  \dot{v} \neq 0
        \label{eq:contact_law_friction}
    \end{gather}
    \noindent
    where $f_{t,u}$ and $f_{t,v}$ are the friction forces along 2 orthogonal directions spanning the plane of contact, $\dot{u}$ and $\dot{v}$ are the corresponding tangential velocities, $\mu$ is the coefficient of friction, $f_{n}$ is the normal contact force. If $\dot{u} = \dot{v} = 0$ and the forces exerted on the element are smaller than the slipping limit, then the friction forces are in equilibrium with those forces and the element is in stick. As for the normal force, it temporally evolves, and in case of separation, there are no friction forces. In this work, we enforce both contact and Coulomb friction using the penalty method \cite{Krack2017}.   
     
    \subsection{Set-up}
    We write the equations of motion of a FEM of a jointed mechanical assembly as
    \begin{equation}
    \mat{M}\ddot{\vec{q}} + \mat{C}\dot{\vec{q}} + \mat{K}\vec{q} + \vt{f}{}{g}(\vec{q}) + \vt{f}{}{J} (\vec{q},\vec{\dot{q}}) = \vt{p}{}{s} + \vt{p}{}{dyn}(\vec{\Omega} t),
    \label{eq:JointEqnsMotion2ndOrder}
    \end{equation}
    \noindent
    where $\vt{f}{}{g} \in \mathbb{R}^{\text{n}}$ consists of geometrically nonlinear forces, $\vt{f}{}{J} \in \mathbb{R}^{\text{n}}$ is the vector of friction and contact nonlinearities, and $\vt{p}{}{s} \in \mathbb{R}^{\text{n}}$ consists of static loads including bolt preclamp forces.     
    
    First, we compute the static configuration $\vt{q}{}{s}$ under the effect of the preclamping forces $\vt{p}{}{s}$. This is given by

    \begin{equation}
        \mat{K}\vt{q}{}{s} + \vt{f}{}{g} (\vt{q}{}{s}) + \vt{f}{}{J}(\vt{q}{}{s}) = \vt{p}{}{s}.      \label{eq:StaticEquilibrium}
    \end{equation} 
    We write the system in $1^{\text{st}}$-order form with the static configurations as the origin by introducing $\tilde{\vec{x}} = \vec{x} - \vt{x}{0}{}$, where $\vec{x} = [\vec{q}, \vec{\dot{q}}]^{\text{T}} \in \mathbb{R}^{2n}$ and $\vt{x}{0}{}=[\vt{q}{s}{}, \vec{0}]^{\text{T}}$. The phase-space representation \eqref{eq:smooth_1st_order} of system \eqref{eq:StaticEquilibrium} is then

   \begin{equation}
        \vec{\dot{\tilde{x}}} =
        \mat{A} [\vec{\tilde{x}} + \vt{x}{0}{}] + \vt{f}{0}{}(\vec{\tilde{x}} + \vt{x}{0}{}) + \vt{f}{s}{} + \epsilon \vt{f}{1}{}(\vec{\Omega}t),
        \label{eq:joint_1st_order}
        \\
        \end{equation}
        \begin{equation*}
        \mat{A} =
        \begin{bmatrix}
            \mat{0} && \mat{I}\\
            -\mt{M}{}{-1}[\mat{K} + \frac{\partial \vt{f}{}{g}}{\partial \vt{\tilde{x}}{}{}}|_{\vec{\tilde{x}}=\vt{0}{}{}} + \mt{K}{}{J}] &&  -\mt{M}{}{-1}\mat{C}
        \end{bmatrix}, \qquad
        \vt{f}{s}{} = \begin{bmatrix}
            \vec{0} \\
            \mt{M}{}{-1} \vt{p}{}{s}
        \end{bmatrix}
        \qquad
        \epsilon \vt{f}{1}{} =
        \begin{bmatrix}
            \vec{0}\\
            \mt{M}{}{-1} \vt{p}{}{dyn}(\Omega t)
        \end{bmatrix},               
        \end{equation*}
    \noindent
    where $\tilde{\vec{x}} = \vec{0}$ is a fixed point of system \eqref{eq:JointEqnsMotion2ndOrder} for $\vt{p}{}{dyn}(\Omega t) = \vec{0}$, $\vt{f}{0}{}$ is the vector of nonlinear forces, and $\mt{K}{}{J}$ consists of the stiffness contributions of the frictional contact elements. The stiffness contribution of a node in contact is

    \begin{equation*}
    \frac{\partial{\vt{f}{i}{J}}}{\partial{\vt{q}{i}{}}} =
        \begin{bmatrix}
            \st{k}{p}{} && 0 && 0 \\
            0 && \st{k}{p}{} && 0 \\
            0 && 0 && \st{k}{p}{}
        \end{bmatrix},
    \end{equation*}
    where $\st{k}{p}{}$ is the numerical penalty parameter used. For contact pairs exhibiting separation at $\tilde{\vec{x}}=\vec{0}$, the element has no contribution to the stiffness matrix. In practice, we perform the spectral analysis of $\mat{A}$ by solving the eigenvalue problem associated with structural dynamics in the configuration space \cite{Thurnher2024}. The spectral analysis of $\mat{A}$ is important for two reasons. First, the linear spectrum informs us about the SSM structure emanating from the fixed point $\vec{\tilde{x}} = \vec{0}$ of the phase space. This is due to the imaginary parts of the eigenvalues of $\mat{A}$ being good approximations of the frequencies of jointed structure at low amplitudes of excitation. Consequently, we can distinguish the slowest SSMs and uncover potential internal resonances. Second, we use the corresponding eigenvectors of $\mat{A}$ to identify the tangent subspace $\mat{V}_{E^{m}}$ at the fixed point for our ROM. We discuss all this in more detail in the next subsection.

    %We note that the existence of the attracting and persistent 2m-dimensional submanifolds in this setting is expected. On one hand, Bettini et al. \cite{Bettini2023} formulated the existence of attracting sets for non-smooth systems. Those are composed of the union of low-dimensional SSMs with switching surfaces. On the other hand, ways to model contact and friction include smoothening the piece-wise forces, in which case the theory of SSMs would be directly applicable. 
    
    \subsection{Data-Assisted Model Reduction using FastSSM and SSMLearn}  
    The mathematical theory of SSMs reviewed in Section \ref{subsection:NonlinearDynamics} is applicable to smooth dynamical systems, with a piecewise smooth extension discussed in \cite{Bettini2023}. In the present work, we will use a data-driven approach in which we find a closely fitting smooth ROM that faithfully reproduces the dominant features of the flow map of the nonsmooth system \eqref{eq:joint_1st_order} in its phase space. As we shall see, this smooth ROM is also predictive: it can predict the forced response of the full nonsmooth system even though it is only trained on trajectory data from the unforced system. 
    
\begin{figure}[H]
        \begin{centering}
        \includegraphics[width=0.5\textwidth]{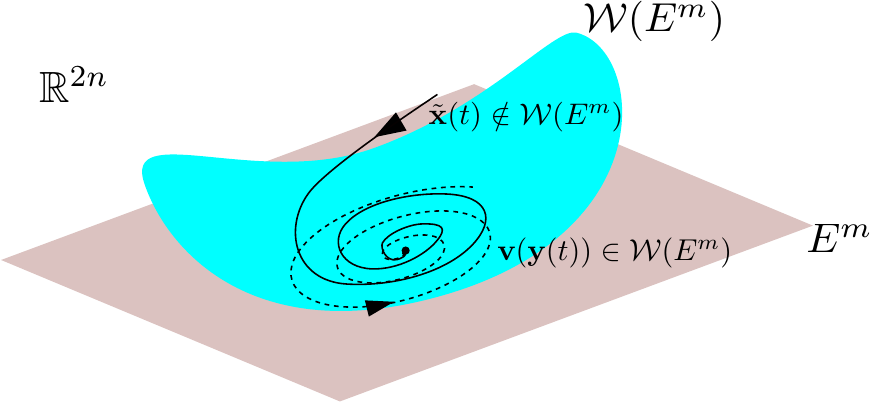}
        \par\end{centering}
        \caption{Sketch of the geometry and dynamics of full and reduced coordinates. A sample trajectory $\vec{\tilde{x}}(t)$ of the full system is sketched in solid line. Its prediction using the SSM-based ROM lies on the $\mathcal{W}(E^{m})$ (in blue) and is sketched in dashed lines.}
        \label{fig:ROM_sketch}
    \end{figure}

    As in the work by Cenedese et al. \cite{Cenedese2023}, we identify an SSM-based ROM using a few transient FE simulations as training data. The initial conditions for these simulations should be selected to activate the dynamic regime of interest. For instance, if we are interested in the first nonlinear modal response of a structure, we then define the initial conditions to be a modulation of the first linear mode. The data that we use for training are the displacements and velocities $\mat{\tilde{X}}\in \mathbb{R}^{2n \times N}$, where $N$ is the total number of time snapshots obtained from the trajectories. Next, we use the open-source algorithms FastSSM \cite{Axas2022} SSMLearn \cite{Cenedese2023} to construct our data-based ROM. This consists of learning the graph of the SSM over the dominant spectral subspace and then learning the dynamics within the SSM. Figure \ref{fig:ROM_sketch} shows a sketch of the relation between the full dynamics, and the SSM obtained in this fashion.\\

    \noindent
    \emph{SSM Geometry from FastSSM}\\
    We compute the reduced coordinates as
    \begin{equation}
           \vec{y} = \vec{w}(\vec{\tilde{x}}) = \mat{V}_{E^{m}}^{\dagger} \vec{\tilde{x}}, 
    \end{equation}
    where $\vec{y} \in \mathbb{R}^{m}$, the columns of $\mat{V}_{E^{m}}$ are the eigenvectors of the slowest m modes of the system, and $(\cdot )^{\dagger}$ denotes the Moore-Penrose pseudo-inverse. As reviewed in Sec. \ref{section:introduction}, based on the theory of SSMs, we anticipate that the targeted dimensionality of the manifold is $m$. Next, a parameterization of the manifold serves to map from the reduced coordinates to the full state-space coordinates. This parameterization is known to admit a Taylor expansion of the form
    \begin{equation}
        \vec{\tilde{x}} \approx \vec{v} (\vec{y}) = \mat{W}\vec{y}^{1:p} = \mat{V}_{E^{m}}\vec{y} + \mat{W}_{2:p} \vec{y}^{2:p},
        \label{eq:param}
    \end{equation}
    \begin{equation*}
        \mat{W} =[\mat{V}_{E^{m}}, \mat{W}_{2},..., \mat{W}_{m}],            
    \end{equation*}
    where $\mt{W}{i}{} \in \mathbb{R}^{n \cross d_{i}}$, with $d_{i}$ being the number of $m$-variate monomials at order $i$, and $p$ is the degree of the polynomial approximation selected for the SSM geometry. The superscript $(\cdot)^{l:r}$ denotes a vector of all monomials at orders $l$ through $r$. For instance, if $\vec{y} = [y_{1}\ y_{2}]^{\text{T}}$, then
    \begin{equation*}
        \vec{y}^{2:3} = \begin{bmatrix}       y_{1}^{2}\ \  y_{1} y_{2}\ \ y_{2}^{2}\ \ y_{1}^{3}\ \ y_{1}^{2} y_{2}\ \ y_{1} y_{2}^{2}\ \ y_{2}^{3} 
        \end{bmatrix}^{\text{T}}.
    \end{equation*}
    \noindent
    To learn this parameterization $\mat{W}$, we use FastSSM \cite{Axas2022} which performs a polynomial regression on the data,
    \begin{equation}
    \mat{W} = \mat{\tilde{X}}(\mat{Y}^{1:p})^{\dagger},
      \end{equation}
    where $\mt{Y}{}{i} \in \mathbb{R}^{m_{i} \cross N}$.\\

    \noindent
    \emph{SSM dynamics from SSMLearn}\\
    Once the SSM geometry is known, we use SSMLearn to find the normal form of the SSM-reduced dynamics \cite{Cenedese2022a}. This form is sparse by construction and is physically interpretable in polar coordinates. Note that this sparsity structure is known and imposed by SSMLearn once the eigenvalues of the linear, autonomous system in \eqref{eq:joint_1st_order} are known. The transformation mapping from normal form coordinates $\vec{z}\in \mathbb{C}^{m}$ to the reduced coordinates $\vec{y} \in \mathbb{R}^{m}$ is denoted by $\vec{y}=\vt{t}{}{}(\vec{z})$. SSMLearn uses a polynomial fit to compute the inverse transformation, which is then given by
    \begin{equation}
        \vec{z} = \vt{t}{}{-1}(\vec{y}) = \mat{H}\vt{y}{}{1:f} = \mt{B}{}{-1}\vec{y} + \mat{H}_{2:f} \vec{y}^{2:f},
    \end{equation}
    where $f$ is the polynomial degree selected for the normal form, which generates the simplified reduced dynamics
    \begin{equation}
        \dot{\vec{z}} = \vec{n}(\vec{z}) = \mat{N}\vec{z}^{1:f} = \mat{\Lambda}\vec{z} + \mat{N}_{2:f} \vec{z}^{2:f}.
        \label{eq:NFmap}
    \end{equation}
    Note that $\mat{N}$ is a sparse matrix that contains nonzero entries only for near-resonant terms \cite{Cenedese2022a}. To compute $\mat{H}$ and $\mat{N}$, SSMLearn solves the minimization problem \cite{SSMLearn}, 
    \begin{equation}
    \underset{\mat{N},\mat{H}}{\text{argmin}} \sum_{j} \left\lVert \frac{\text{d}}{\text{d}\vec{y}}\vt{t}{}{-1}(\vt{y}{j}{})\dot{\vec{y}}_{\text{j}}-\vec{n}(\vt{t}{}{-1}(\vt{y}{j}{}))\right\rVert^{2},
    \end{equation}
    where $\vec{y}_{j} \in \mathbb{R}^{m\times L}$, where $L$ is the number of training data points.
    
    The SSM-based ROM predicts the dominant dynamics of all unforced trajectories with generic initial conditions close to the SSM. Importantly, this ROM can also be used to predict the forced response of the full system \eqref{eq:joint_1st_order}. Specifically, if the time-dependent forcing is quasi-periodic as in \eqref{eq:joint_1st_order}, forcing can be added at leading order to the unforced ROM \eqref{eq:NFmap} as \cite{Breunung2018a} 
    \begin{equation}
        \dot{\vec{z}} = \vec{n}(\vec{z}) + \mt{B}{}{-1} \mat{V}_{E^{m}}^{\dagger} \vt{p}{1}{dyn} (\vec{\Omega} \text{t}).
        \label{eq:NF_forced}
    \end{equation}
    In this work, we use Eq. \eqref{eq:NF_forced} and use the autonomous parameterization of the SSM to predict the forced response of the full FEM. While not necessary here, the accuracy of the forced response prediction can be increased by using the full time-dependent parameterization of the SSM instead of Eq. \eqref{eq:param}.\\

    \noindent
    \emph{Error Metric}\\
    To compute the error of the ROM \eqref{eq:NFmap}, we use the Normalized Mean-Trajectory Error (NMTE) error \cite{Cenedese2022a,Cenedese2022}, defined as

    \begin{equation}
        \text{NMTE} = 100\frac{ \sum^{L}_{k=1} \lVert \vec{\tilde{x}}(t_{k}))-\vec{\tilde{x}}_{\text{rec}}(t_{k}))\rVert}{L||\vt{\underline{x}}{}{}||},
        \label{eq:NMTE}
    \end{equation}
    where $L$ is the number of observations made along the trajectory $\vec{\tilde{x}}(t)$, $\vec{\tilde{x}}_{\text{rec}}(t)$ is the reconstructed trajectory, and $\vt{\underline{x}}{}{}$ is a normalization vector. In this work, $\vt{\underline{x}}{}{}$ is the trajectory vector at the time sample where the maximum norm is reached, i.e. $\underset{k}{\text{max}}||\vec{\tilde{x}}(t_{k})||$. We judge the ROM to be accurate if the NMTE error less than 10\% on the set of trajectories used in the validation (but not in the training) \cite{Cenedese2023}.
        \section{TRCBenchmark Test Case}
    \label{section:TRCBenchmark}
\begin{figure}[H]
        \centering
        \begin{minipage}{0.4\textwidth}
            \centering       \includegraphics[scale=0.25]{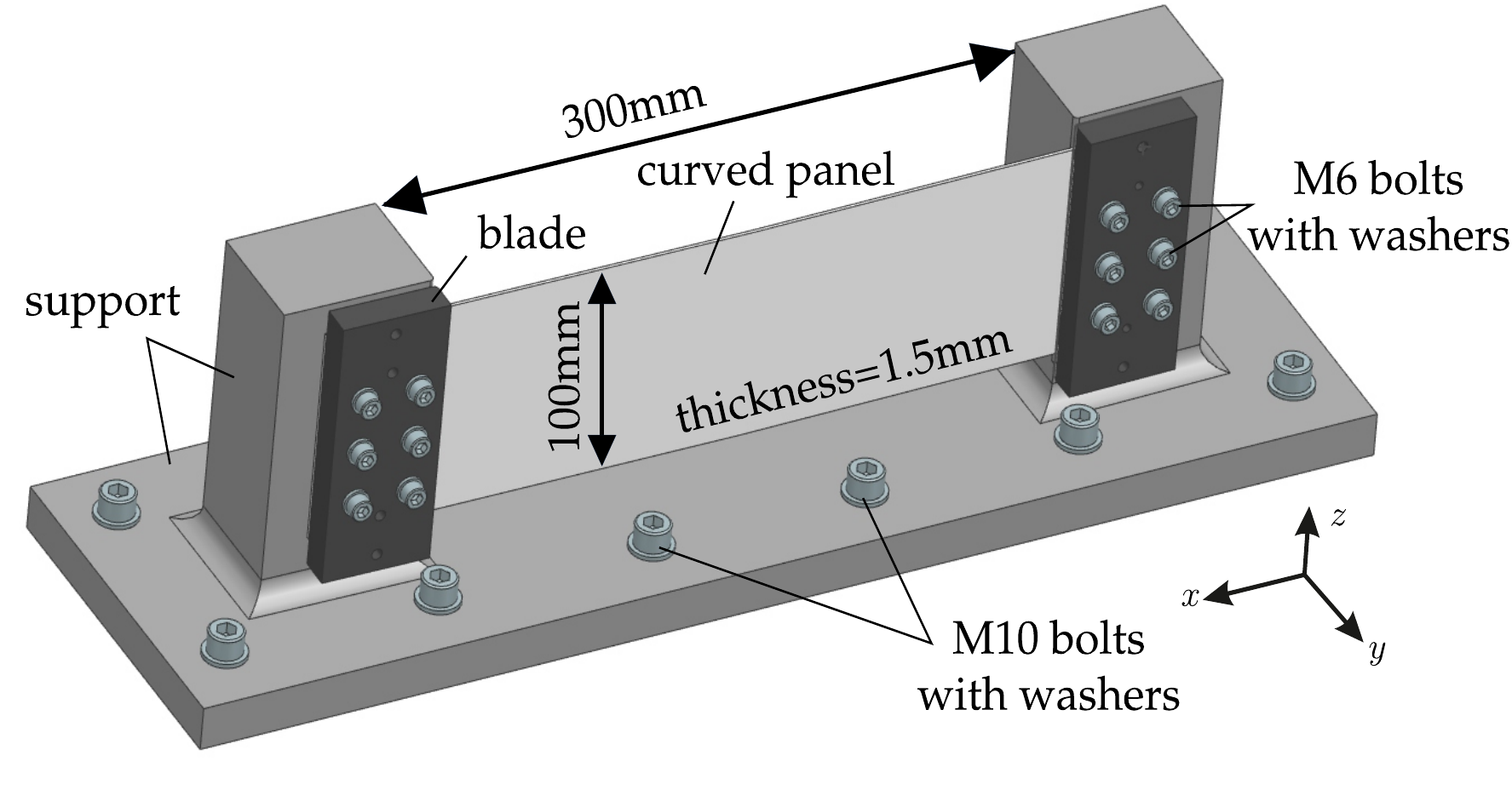}
    \caption{CAD model of the TRC Benchmark Structure \cite{darus-3147_2022}.}
    \label{fig:CAD}
        \end{minipage}\hfill
        \begin{minipage}{0.4\textwidth}
            \centering
           \includegraphics[scale=0.4]{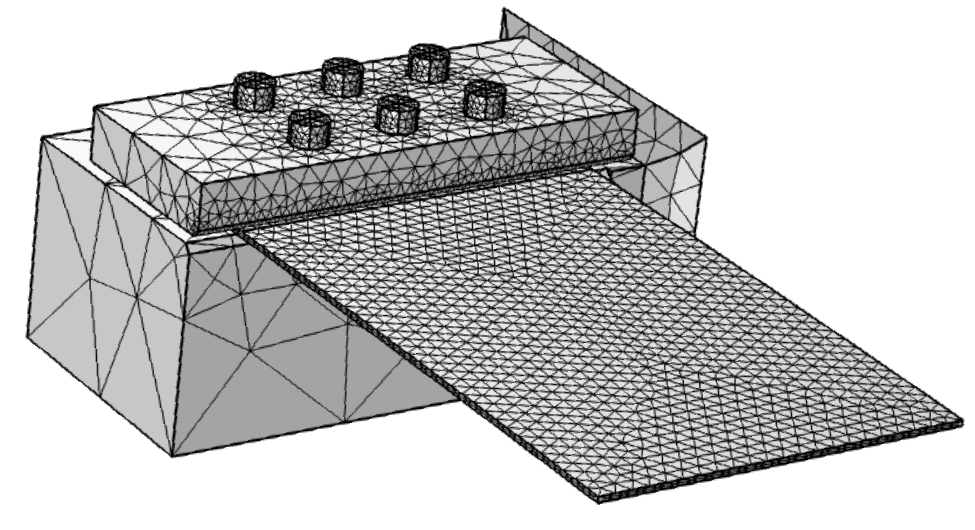}
    \caption{Finite element model. Number of degrees of freedom = 187,920.  Friction is modelled as Coulomb friction with $\mu=0.6$.}
    \label{fig:FEM}
        \end{minipage}
    \end{figure}

    \emph{Model}\\
    Figure \ref{fig:CAD} shows a CAD model of the TRC Benchmark structure \cite{darus-3147_2022}. It consists of a slender panel, a support structure featuring two vertical cantilevers, and two steel blades. Six bolts pre-clamp the components together as shown in the figure. The panel is made of stainless steel 301-1.4310, while the blade and the support structure are made of hardened steel 1.7147, and the bolts are ISO 4762 - M6. The structure is subjected to a base excitation in the $y$-direction, with a frequency triggering the first nonlinear modal response at growing amplitudes. We define the geometry of the structure by importing the CAD model \cite{darus-3147_2022} to COMSOL Multiphysics$\textsuperscript{\textregistered}$ 6.1 to construct the FEM and run the simulations.

    Figure \ref{fig:FEM} shows our FEM of the half-system, where we have imposed symmetry boundary conditions at the middle of the plate. The mesh contains 37,316 quadratic tetrahedral elements with a total of 187,920 degrees of freedom. Table \ref{table:materials} shows the material properties assigned to the different components. We use the Green-Lagrange strain tensor to model the geometric nonlinearity of the structure. As for modeling the joint, we define frictional contact at the following interfaces: between the support and panel, between the panel and the blade, and between the panel and the bolt heads. Friction is modeled as Coulomb friction with an assumed coefficient of friction $\mu$ = 0.6 {\cite{Schwingshackl2012}. Both friction and contact are enforced using the penalty method \cite{Krack2017}. We impose an initial stress of around $416.5\cross10^{6} \frac{\text{N}}{\text{m}^2}$ in the shafts of each bolt to generate a bolt preclamping force of 9.1 KN. This force is computed via the empirical formula \cite{bolt_formula}
    
    \begin{equation}
        F = \frac{T}{(0.159P + 0.578d_{2}\mu_{T} + 0.5D_{f}\mu_{H})},
        \label{eq:bolt_force}
    \end{equation}
    \noindent
    where the applied torque T = 10.1 Nm\cite{darus-3147_2022}, the bolt pitch P = 1 mm, the nominal bolt diameter $d_{2}$ = 6 mm, the average of the bolt cap and shaft diameters $D_{f}$ = 8 mm, the bolt head friction coefficient is assumed $\mu_{H}$ = 0.2 \cite{Wall2022}, and the thread friction coefficient $\mu_{T}$ = 0.2. We use linear Rayleigh damping with a modal damping ratio of 0.4$\%$ for the first two modes.
    
\begin{table}[H]
\centering
\begin{tabular}{|c|c|c|c|}
\hline                           & \multicolumn{1}{l|}{$\rho [\frac{Kg}{m^3}]$} & \multicolumn{1}{l|}{E [GPa]} & \multicolumn{1}{l|}{$\nu$} \\ \hline
Panel                             & 7,900          & 195                              & 0.29                                   \\ \hline
Blade & 7,900                                    & 200                              & 0.29                           \\ \hline
\multicolumn{1}{|l|}{Cantilever Support} & 7,900   & 200                              & 0.29                                    \\ \hline
Bolts                                    & 7,850                 & 210                              & 0.3                      \\ \hline
\end{tabular}
\caption{Material properties of the FEM.}
\label{table:materials}
\end{table}
    
%The panel is has $\text{E}=195$ GPa, $\nu=0.29$, and $\rho=7900\ \frac{\text{Kg}}{\st{m}{}{3}}$. The blade has $\text{E}=200$ GPa, $\nu=0.29$, and $\rho=7900\ \frac{\text{Kg}}{\st{m}{}{3}}$. The bolts have $\text{E}=210$ GPa, $\nu=0.3$, and $\rho=7850\ \frac{\text{Kg}}{\st{m}{}{3}}$.  The analysis is nonlinear due to both geometric and frictional contact nonlinearity.  Geometric nonlinearity is present due to the slenderness of the panel. Frictional contact is present at the interfaces joining the different components.

    \noindent
    \emph{Data and ROM construction}\\
    After computing the static equilibrium of the structure under the effect of bolt pre-clamp forces, we extract the system matrices from COMSOL and perform the spectral analysis of the matrix $\mat{A}$ defined in Eq. \eqref{eq:joint_1st_order}. The first two linear modes of the structure are a bending mode of the panel at 109.4 Hz and a torsional mode at 196.6 Hz. The corresponding mode shapes computed from the undamped eigenvalue problem are plotted in Figs. \ref{fig:BendingTime}-b. Since we aim to make predictions of the forced dynamics triggering the first mode, we a priori envision constructing a two-dimensional SSM serving as a nonlinear continuation of the corresponding 2D linear spectral subspace of the linearized system at the origin. To acquire the training data needed, we run a freely decaying FE simulation with initial conditions perturbations of the first linear mode shape. The maximum displacement using this perturbation is chosen to activate the nonlinear dynamics of the system up to a displacement of 2mm. Figure \ref{fig:BendingTimeRef} shows the time response of the corner point P, while we demonstrate the frequency content of this signal using short-time fourier transform in Fig. \ref{fig:BendingFreqRef}. We note that, in addition to the higher harmonics of the $1^{\text{st}}$ linear frequency, the frequency of the $2^{\text{nd}}$ mode is also present, indicating a 1:2 internal resonance. This can be seen from the evolution of the thick frequency band at around 200 Hz into two distinct frequency components as the response decays. The existence conditions for SSMs \cite{Haller2016} require that we select the 4D spectral subspace containing this resonance and construct a 4D SSM tangent to this subspace at the origin of system \eqref{eq:joint_1st_order}.

     \begin{figure}
        \centering
        \begin{subfigure}{0.5\textwidth}
            \centering
            \includegraphics[width=1\textwidth]{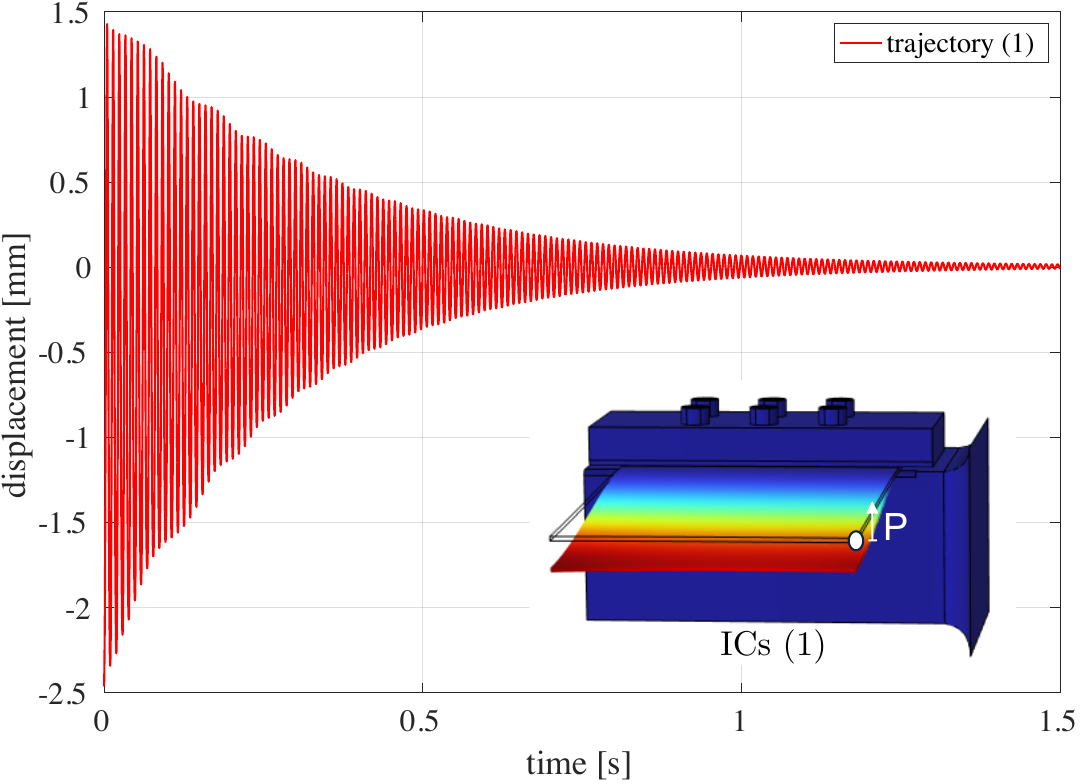}
            \caption{}
            \label{fig:BendingTimeRef}
        \end{subfigure}\hfill
        \begin{subfigure}{0.5\textwidth}
            \centering
            \includegraphics[width=1\textwidth]{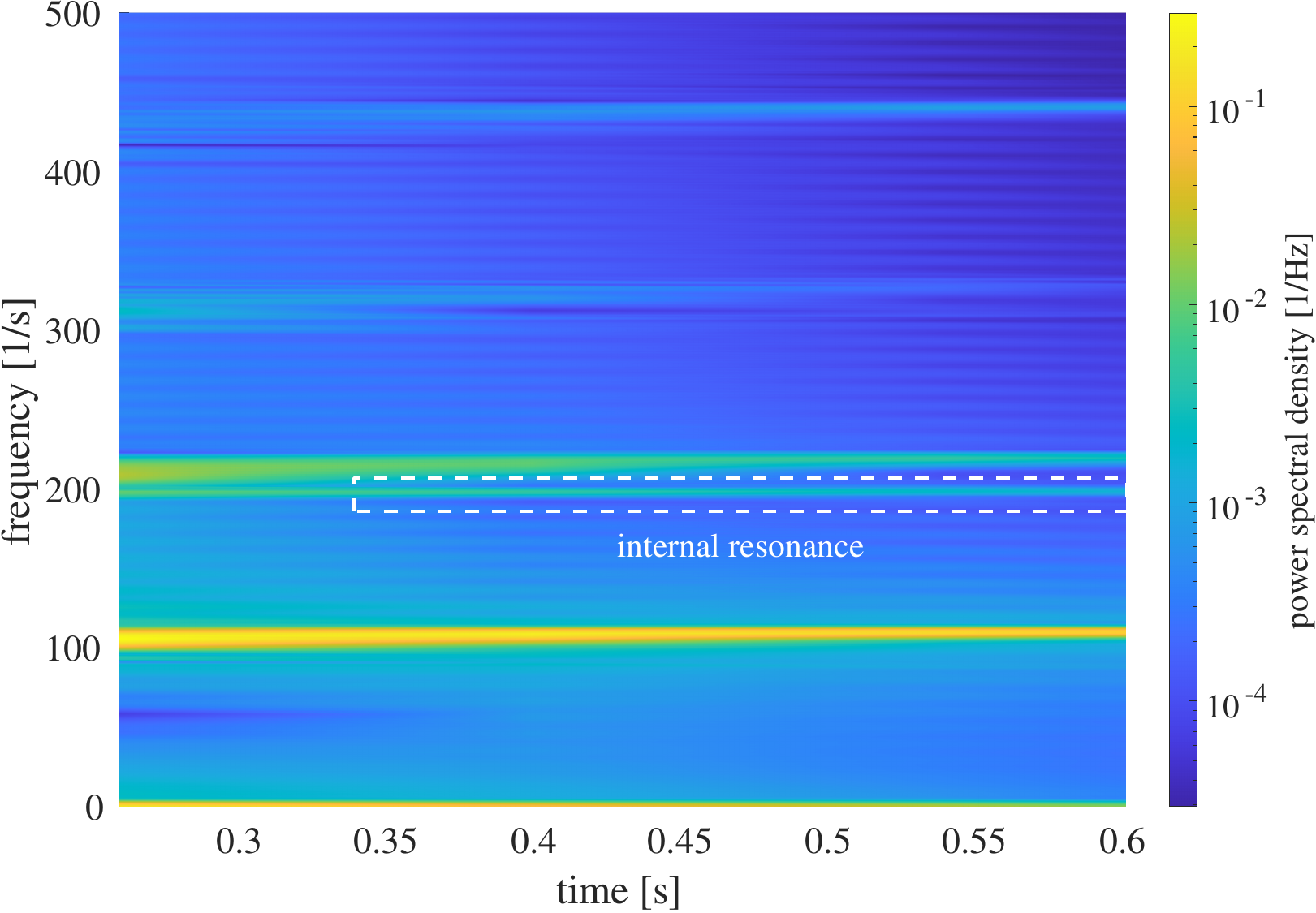}
            \caption{}  \label{fig:BendingFreqRef}
        \end{subfigure}
        \caption{Plot (a) shows the transient time response of trajectory (1) at DOF P. Plot (b) shows its spectrogram. We note an internal resonance of the structure at high amplitudes of vibration as the second harmonic component of the first bending mode activates the torsional mode. We therefore target a 4-dimensional SSM.}
    \end{figure}

             \begin{figure}
        \begin{subfigure}{0.5\textwidth}
            \centering
            \includegraphics[width=1\textwidth]{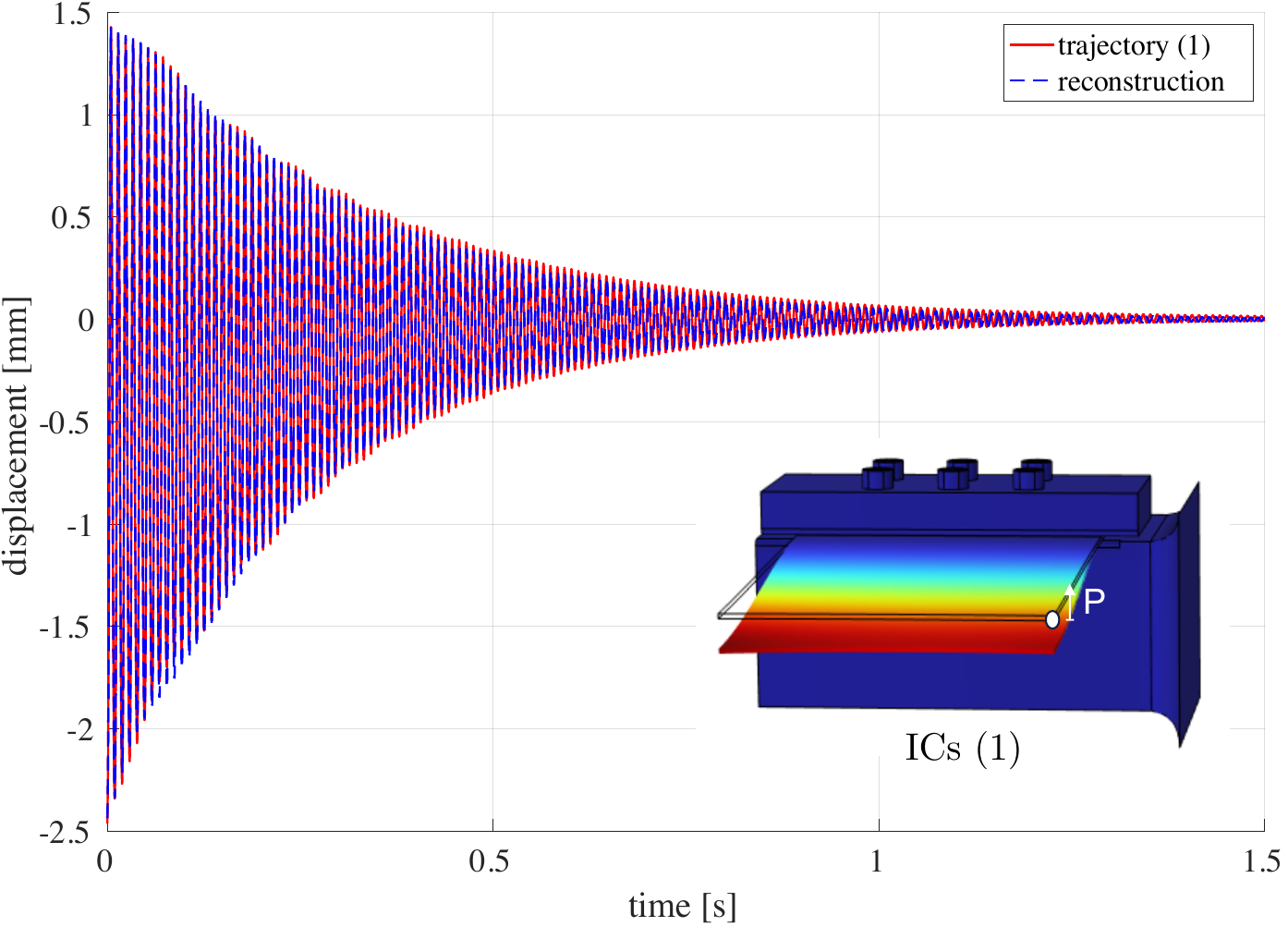}
            \caption{}
            \label{fig:BendingTime}
        \end{subfigure}\hfill
        \begin{subfigure}{0.5\textwidth}
            \centering
            \includegraphics[width=1\textwidth]{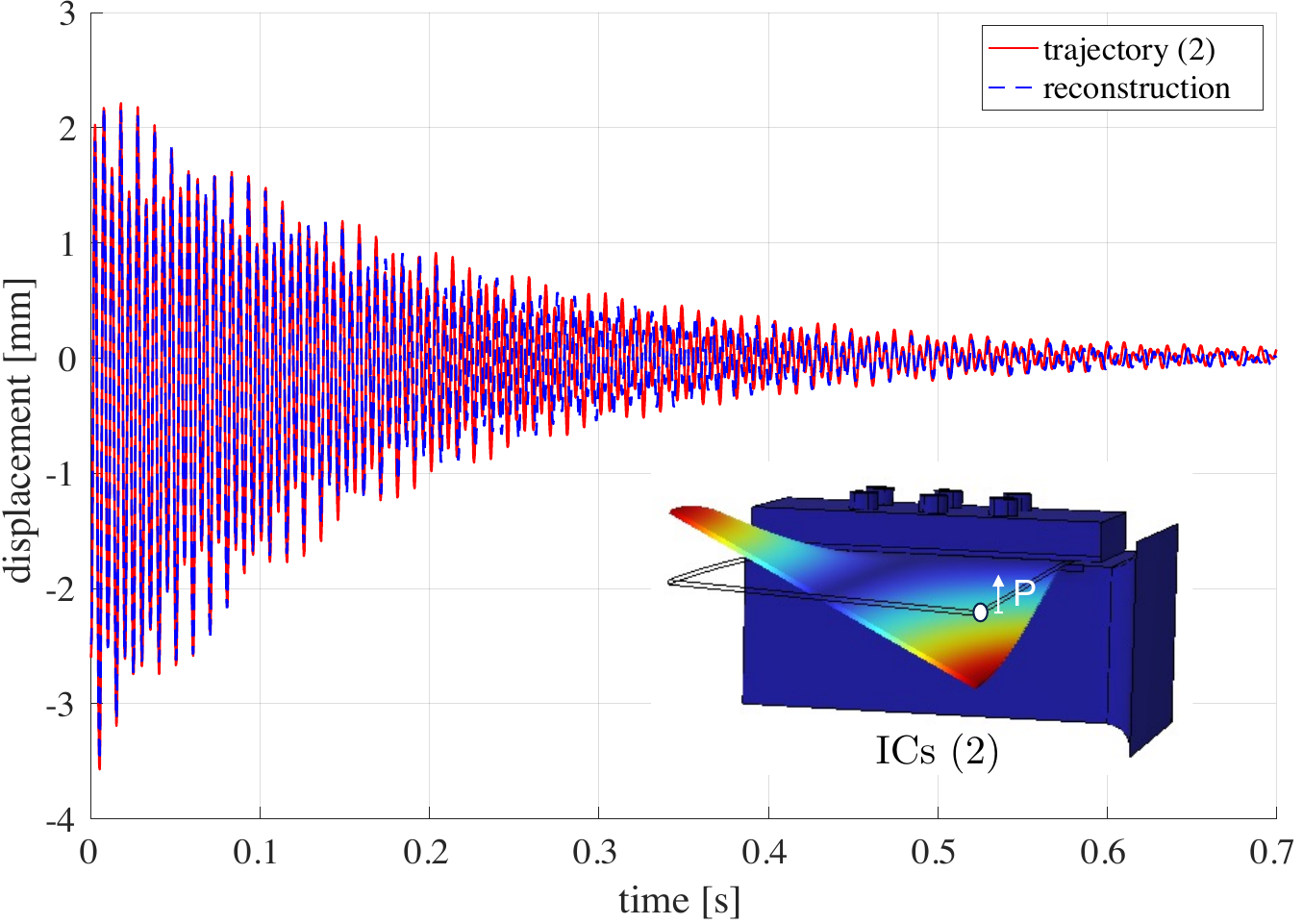}
            \caption{} 
            \label{fig:TorsionTime}
        \end{subfigure}
        \vfill         
        \centering
        \begin{subfigure}{0.5\textwidth}
            \centering
         \includegraphics[width=1\textwidth]{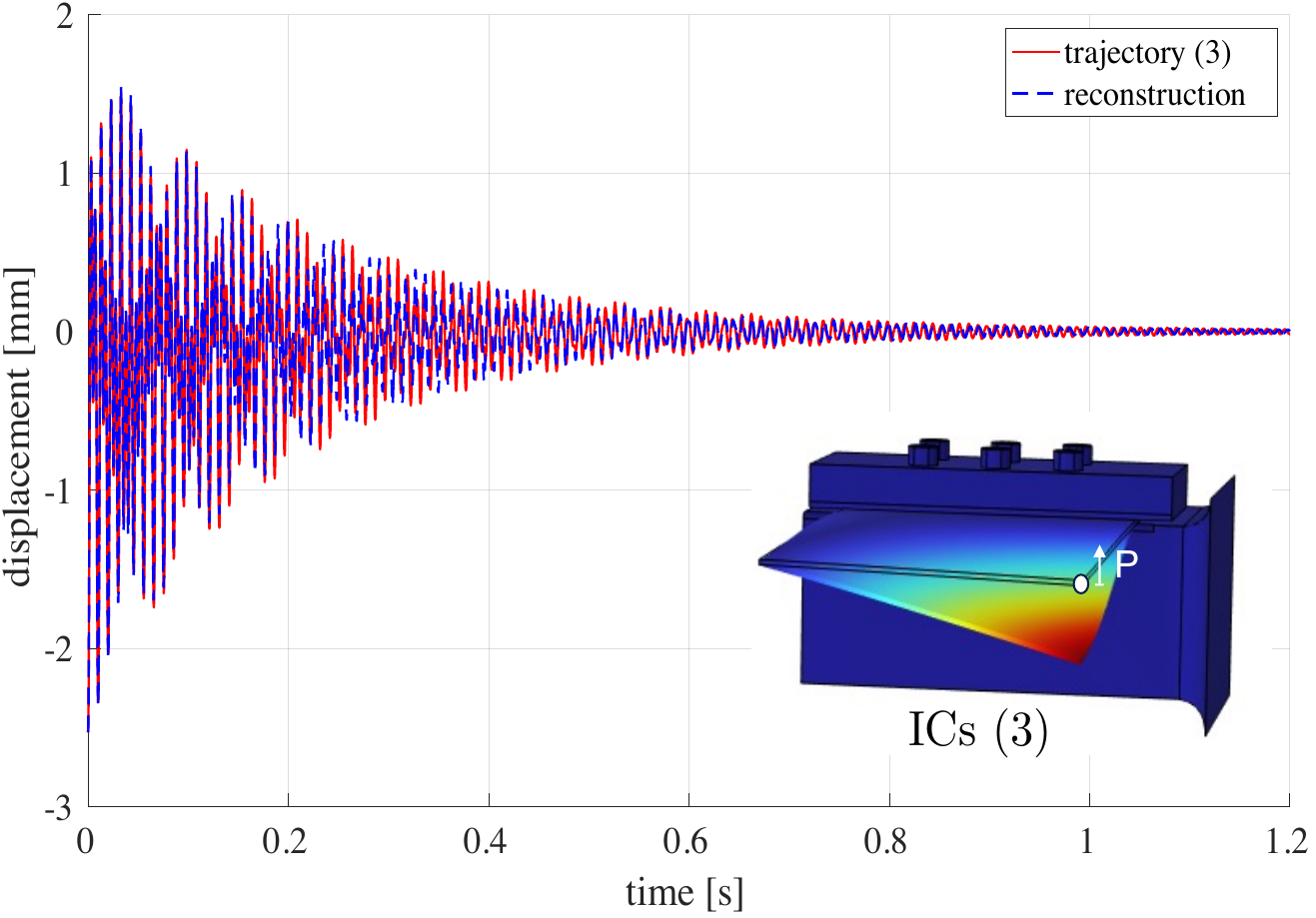}
        \caption{}
            \label{fig:MixTrainTime}
        \end{subfigure}\hfill
        \begin{subfigure}{0.5\textwidth}
            \centering
            \includegraphics[width=1\textwidth]{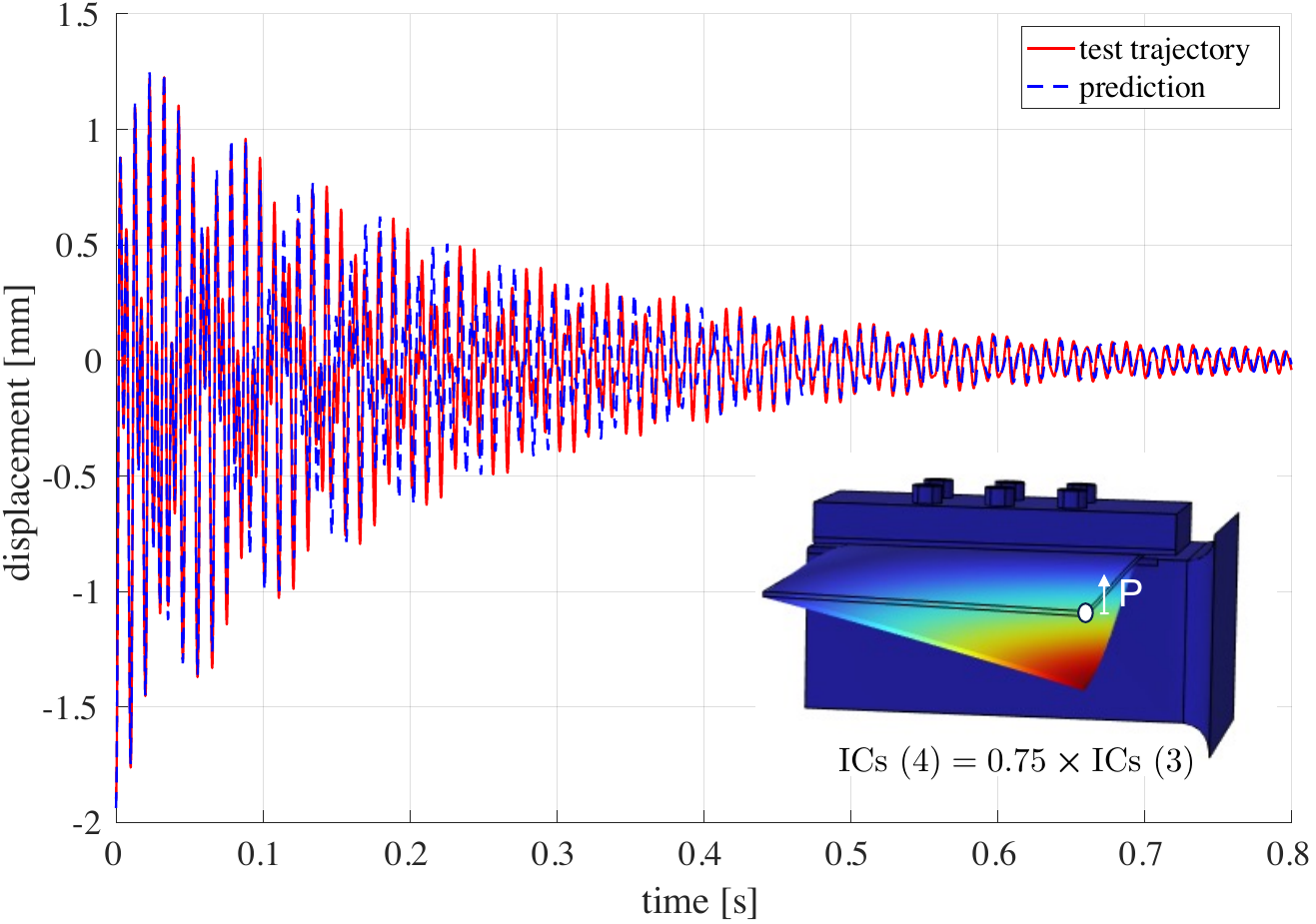}
            \caption{}
          \label{fig:MixTestTime}
        \end{subfigure}
        \caption{Plots (a), (b), and (c) show the free-decaying response of the point P for trajectories (1), (2), and (3) and their reconstruction using the 4D SSM-reduced ROM. Plot (d) shows the prediction of the time response of the test trajectory not used in the training of the SSM-reduced model.}
    \end{figure}    

         \begin{figure}
        \centering
        \begin{subfigure}{0.5\textwidth}
            \centering
         \includegraphics[width=1\textwidth]{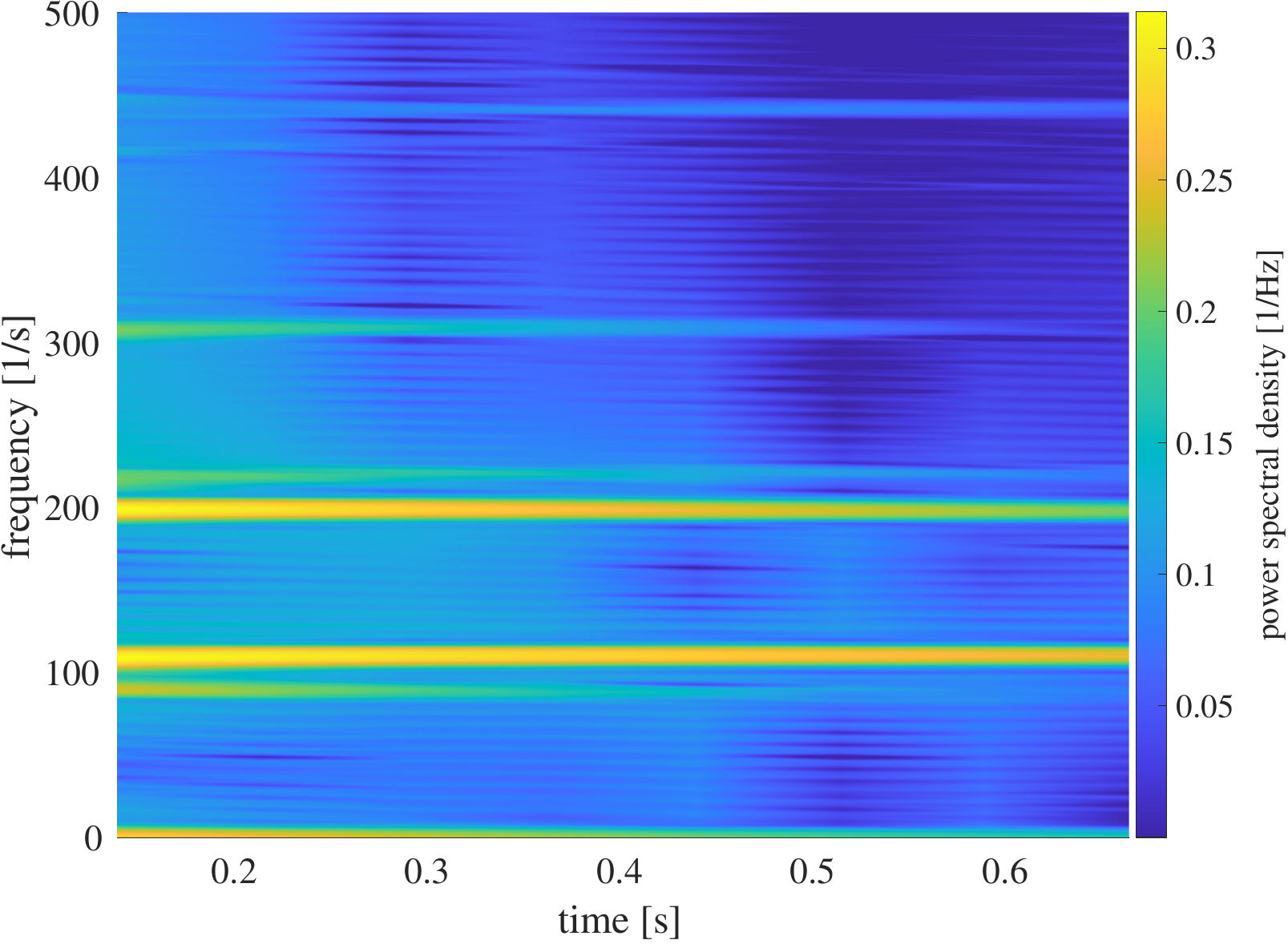}
            \caption{}
            \label{fig:TestFreqReference}
        \end{subfigure}\hfill
        \begin{subfigure}{0.5\textwidth}
            \centering
            \includegraphics[width=1\textwidth]{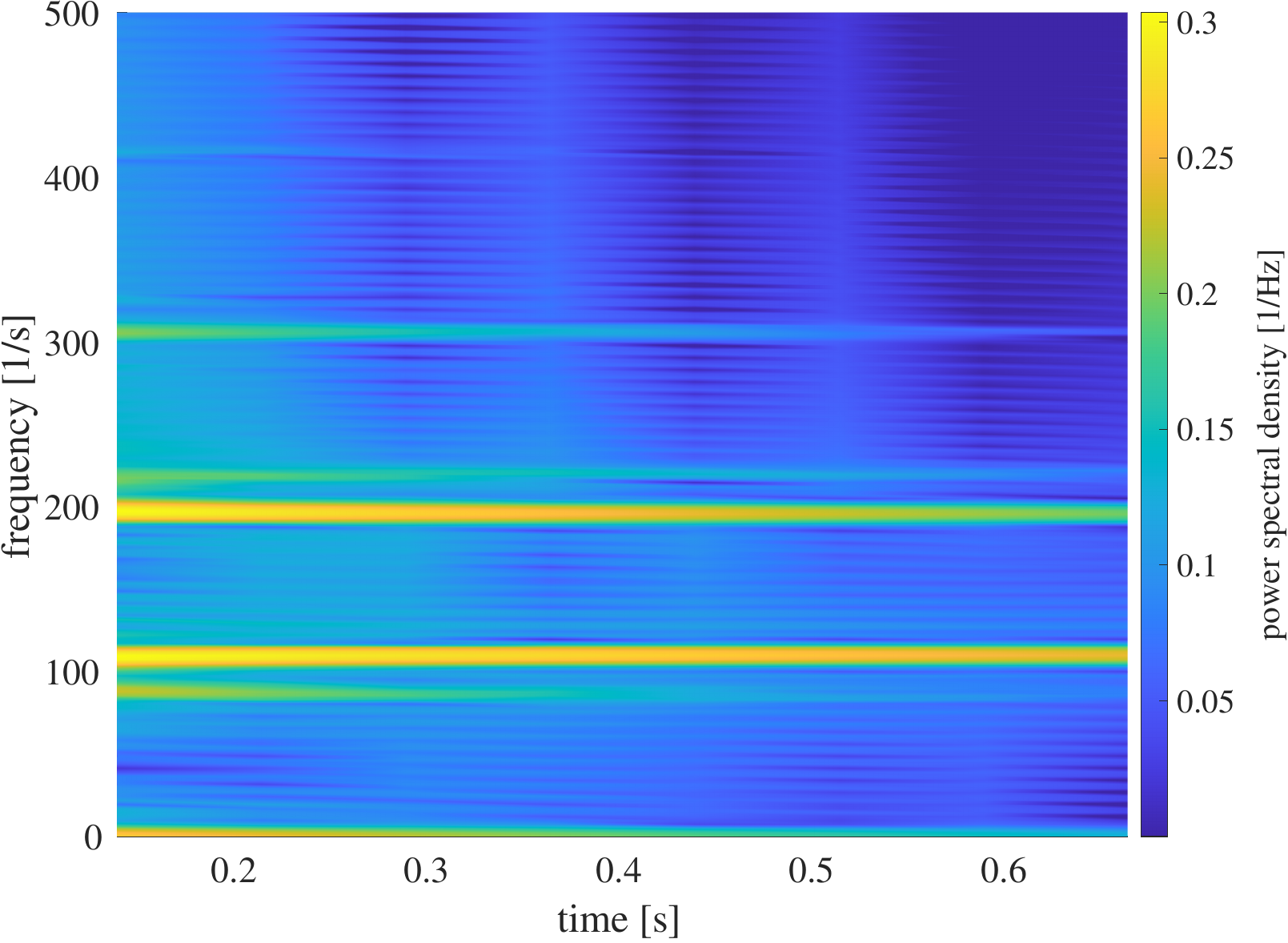}
            \caption{}  \label{fig:TestFreqPrediction}
        \end{subfigure}
        \caption{Plot (a) shows the spectrogram of the free-decaying evolution of the corner point P shown in Fig.d. Plot (b) shows the prediction of this evolution using the SSM-based ROM. The ROM successfully predicts the frequency content of the vibrations, including the $3^{\text{rd}}$ harmonic component of the bending mode at around 310 Hz, and a subharmonic component of the torsion mode below 100 Hz.}
    \end{figure}    
    
    Accordingly, we launch two additional freely decaying trajectories to generate data that covers the neighborhood of the 4-dimensional SSM constructed over the two dominant modes of the linearized system. The initial conditions of the first additional trajectory are perturbations of the torsional mode. The second additional training trajectory has initial conditions equal to a linear combination of the first two modes. Figures \ref{fig:TorsionTime} and \ref{fig:MixTrainTime} show the time evolution of the corresponding trajectories in red. We use the displacements and velocities of the three training trajectories as input to FastSSM and SSMLearn to construct a 4D, SSM-based ROM, as reviewed in Section \ref{section:joints}. Based on the NMTE error \eqref{eq:NMTE}, we select polynomials of order 7 for the geometry and its reduced dynamics, which results in a 9.16$\%$ NMTE error. The corresponding reconstructions of the training trajectories are shown in blue in Figs. \ref{fig:BendingTime}-c. Below we write out the 4D SSM-reduced model in polar coordinates $(\rho_{j},\theta_{j})$ defined via $z_{j}=\rho_{j}e^{i\theta_{j}}$, up to second order, while the complete normal form up to order six is given in Appendix A1:

\begin{minipage}{0.45\textwidth}
  \begin{align*}
\dot{\rho}_{1}\rho^{-1}_{1} = &-3.51-8.79 \rho^{2}_{1}-75.80 \rho^{2}_{2}\\  &+\mathrm{Re}((-19.76+27.61i) \rho_{2} e^{i(-2 \theta_{1}+ \theta_{2})})\\& + \mathcal{O}(||\vec{\rho}||^{3}),\\ \dot{\rho}_{2}\rho^{-1}_{2} = &-5.89-19.06 \rho^{2}_{1}-32.94 \rho^{2}_{2}\\  &+\mathrm{Re}((-0.0063+11.62i) \rho^{2}_{1} \rho^{-1}_{2} e^{i(2 \theta_{1}- \theta_{2})})\\& + \mathcal{O}(||\vec{\rho}||^{3}),\\
    \end{align*}
\end{minipage}
\begin{minipage}{0.50\textwidth}
\begin{equation}
\begin{aligned}
\dot{\theta}_{1} = &+692.94-904.69 \rho^{2}_{1}-287.24 \rho^{2}_{2}\\ &+\mathrm{Im}((-19.76+27.61i) \rho_{2} e^{i(-2 \theta_{1}+ \theta_{2})}) \\& + \mathcal{O}(||\vec{\rho}||^{3}), \\ \dot{\theta}_{2} = &+1231.06+742.60 \rho^{2}_{1}+602.42 \rho^{2}_{2}\\ &+\mathrm{Im}((-0.0063+11.62i) \rho^{2}_{1} \rho^{-1}_{2} e^{i(2 \theta_{1}- \theta_{2})})\\& + \mathcal{O}(||\vec{\rho}||^{3}).\\
    \end{aligned}
    \label{eq:ROM_NF}
    \end{equation}
\end{minipage}

    For validation, we run a test free-decaying simulation with initial conditions equal to a different linear combination of the first two modes. As shown in Fig. \ref{fig:MixTestTime}, the SSM-based reduced model accurately predicts the evolution of this decaying test trajectory. Further, the spectrograms of the two responses shows accurate frequency-domain predictions, as shown in Figs. \ref{fig:TestFreqReference} and \ref{fig:TestFreqPrediction}. The simulations were run on iMac Pro 2.3 GHz 18-Core Intel Xeon W. The average time taken for a training free-decay simulation was 22.4 days. Note, however, that the software uses implicit generalized alpha-method for time integration, which does not exploit the parallel computing resources available. The time step size is set as free, and thus is adaptively determined to satisfy the user-set relative tolerance of $1\cross 10^{-5}$ \cite{comsolreference}. These large computational times can be significantly reduced in case an explicit time integration method is used\footnote{An explicit time integration of the FE model was not pursued in this work since the available FE software, COMSOL Multiphysics$\textsuperscript{\textregistered}$ 6.1, does not support explicit time integration for contact problems.}. For instance, Sandia National Labs reported in \cite{TRCarticle} a computational time of around only 50 hours for a ring-down simulation of a high-fidelity FEM of the same structure on a high-performance computing cluster. One must also appreciate that once the training data is generated, the full forced response of the system can be predicted for arbitrary types of forcing. We illustrate this next for periodic forcing.\\
    
    \noindent
    \emph{Predicting the periodically forced response}\\
    Having successfully predicted the nonlinear dynamics of trajectories not used in the training of the SSM-based ROM, we are ready to forecast the response of the full system under periodic external forcing. For periodic base excitation in the $y$-direction, the external forcing vector is
    \begin{equation}
        \vt{p}{}{dyn} \cos(\Omega \text{t})= -\text{a}\mat{M}\vec{b} \cos(\Omega \text{t}),
        \label{eq:base_excitation}
    \end{equation}
    where $a$ is the magnitude of base acceleration, $\mat{M}$ is the mass matrix, $\vec{b}$ is a vector with unit entries for DOFs in the $y$-direction DOFs and zeros elsewhere. The force vector $\vt{p}{}{dyn}$ is projected onto the tangent space and transformed to normal form as in \eqref{eq:NF_forced}. We compute the periodic orbits of the 4D SSM-reduced model \eqref{eq:ROM_NF} using the COCO toolbox \cite{doi:10.1137/1.9781611972573}. Figure \ref{fig:FRFsBackbone} shows the frequency response curves of the structure for different base excitation amplitudes, all together computed in less than 2 minutes. In contrast, given a single frequency and a single amplitude of excitation, we compute the corresponding single periodic response using a full FE simulation in more than 10 days. In total, we performed 6 such full simulations for validation. As shown in Fig. \ref{fig:FRFsVsComsol}, there is a close overall agreement between the forced response of the SSM-based ROM \eqref{eq:ROM_NF} and the full FEM. 
    
    Next, we compare our predictions results with the available experimental forced response data \cite{TRCarticle}. We plot the experimental results from the 25 tests performed on three different, imperfect, specimens of the slender panel. The variability of experimental results is explained in detail in \cite{TRCarticle,Bhattu2024}. In this work, we aim to compare our prediction with the overall trend of the experimental modal characteristics. Therefore, in Fig. \ref{fig:FRFsBackbone}, we compare the normalized amplitude-dependent resonance frequencies. As shown, we correctly predict the softening behavior observed in experiments up to displacement around 1.5 mm. 
    
    In the experiments, the modal damping ratio was obtained by the period-averaged power balance between dissipation and the base excitation inertial forces \cite{Bhattu2024,Muller2022}. Here, we compute the modal damping prediction by evaluating the instantaneous modal damping of the first mode along the resonance curve, in polar normal form coordinates, at the different displacement levels. The full expression of this instantaneous modal damping is given by the right-hand side of \eqref{eq:NF_damping1st} in Appendix A1. Figure \ref{fig:ModalDamping} shows good agreement with the experimental results. Finally, we use the ROM in eq.\eqref{eq:ROM_NF} to predict friction stresses on the interface. We plot time snapshots of contact and friction stresses due to the periodic response encircled in Fig. \ref{fig:stresses}. These plots highlight that the reliability of our predictions extends to local stresses inside the joint as well. Additionally, to illustrate the microslip dynamics of this jointed structure, we show in Appendix A2 the activation of contact-separation nonlinearities as well as frictional slippage on interfaces.
    
    Next, we compare our results with the blind predictions presented in the TRChallenge \cite{TRCarticle}. Figures \ref{fig:BlindFrequency} and \ref{fig:BlindDamping}, originally presented in \cite{TRCarticle}, now additionally include our SSM-based predictions for the backbone curves and damping ratios. As seen in the figures, our predictions lie within the experimental range of the response exhibited by the structure. All modeling and ROM details of the other methods have been presented in \cite{TRCarticle}. While the equation-driven approach by ETH (in pink) appears to predict well the experimental modal damping, the backbone curve is overly stiff, indicating that the dynamics are not correctly captured\footnote{Further investigations identified the cause in the insufficient number of modes retained in the projection subspace.}. In addition, the backbone curve computed via the UW-BYU (blue curve) approach is also within the experimental range, but the prediction of this approach for the modal damping ratios is clearly inconsistent. This approach fits a single-degree-of-freedom ROM, assumes mono-harmonic response, and identifies the geometric and friction nonlinearities separately. In view of this, any particular success this method had on this example is not expected to generalize to other problems. In contrast, we make no assumptions about the response and construct our SSM-based ROM in a dynamics-based, systematic way independently of the specifics of the jointed structure under study here.

             \begin{figure}[H]
        \centering
        \begin{subfigure}{0.5\textwidth}
            \centering
         \includegraphics[width=1\textwidth]{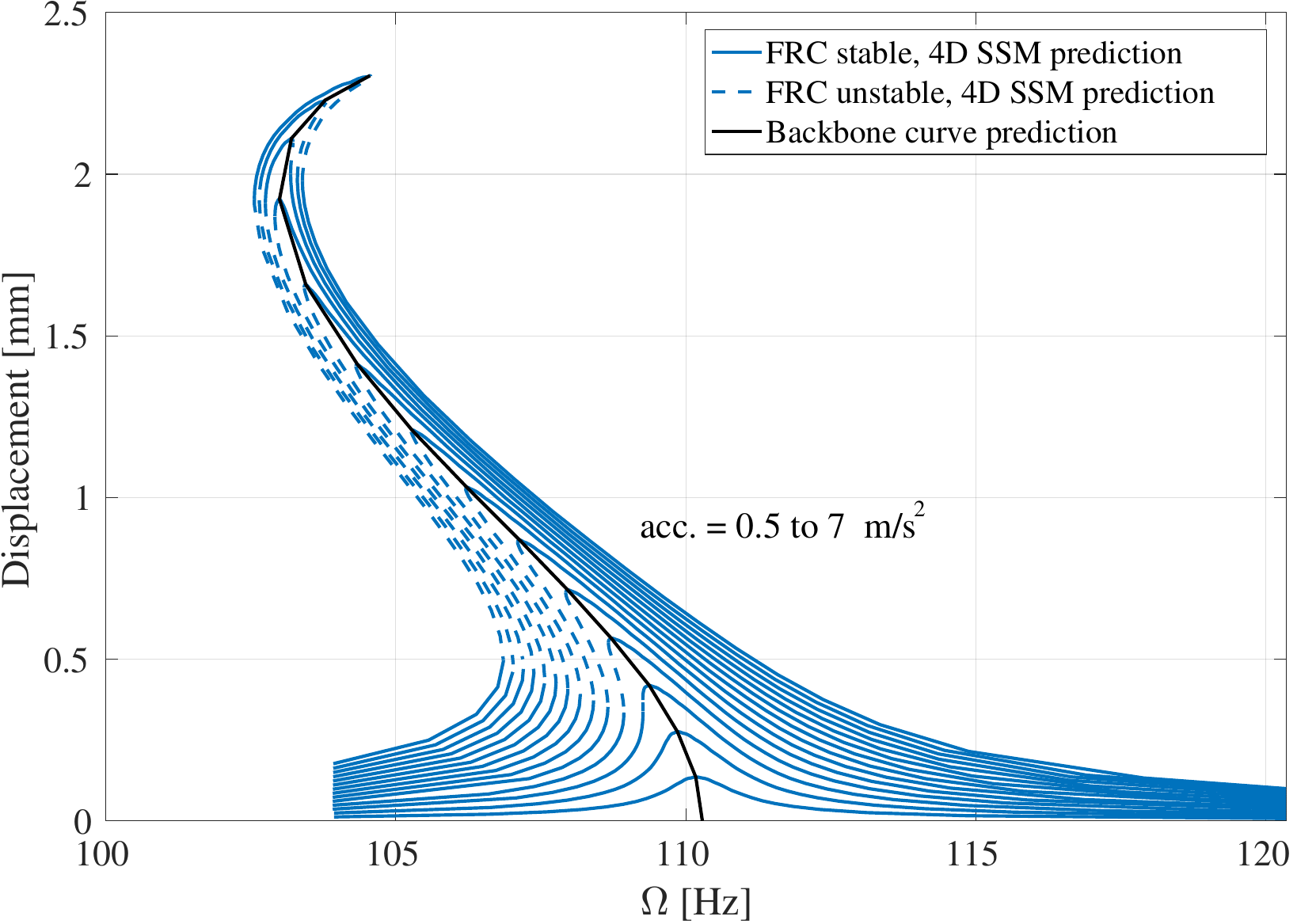}
            \caption{}
            \label{fig:FRFsBackbone}
        \end{subfigure}\hfill
        \begin{subfigure}{0.5\textwidth}
            \centering
            \includegraphics[width=1\textwidth]{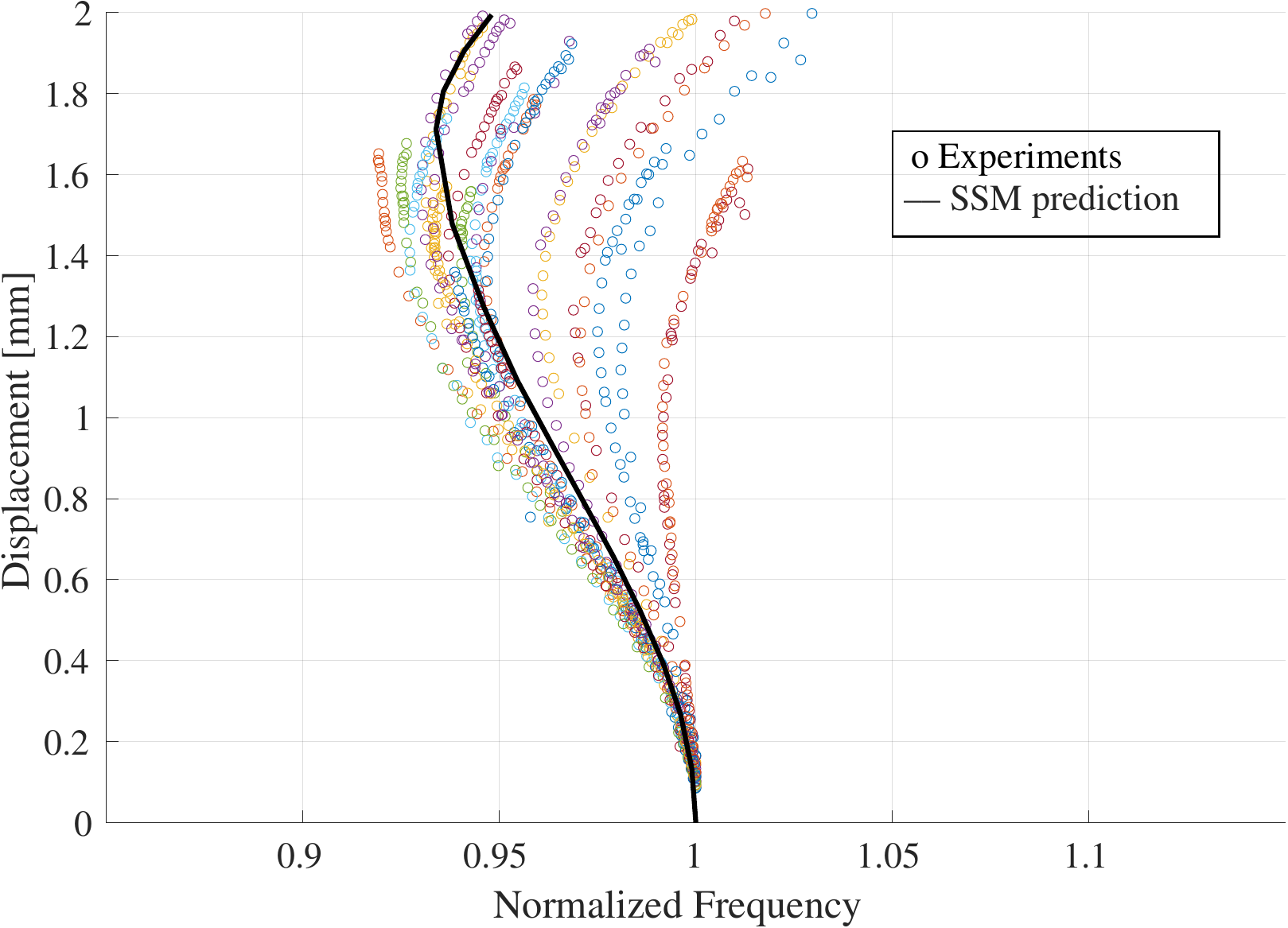}
            \caption{}  \label{fig:AllExperimentsNormalized}
        \end{subfigure}
        \caption{Plot (a) shows SSM-based predictions of maximum out-of-plane displacement responses at the center of the panel due to increasing base acceleration levels. Plot (b) shows the comparison between normalized SSM-predicted resonance curves and experiments. As defined experimentally, the displacement measure here is $\sqrt{2} \cross \text{root mean square of the out-of-plane displacements at the center of the panel across one time period}$.}
    \end{figure}

         \begin{figure}[H]
        \centering
        \begin{subfigure}{0.5\textwidth}
            \centering
         \includegraphics[width=1\textwidth]{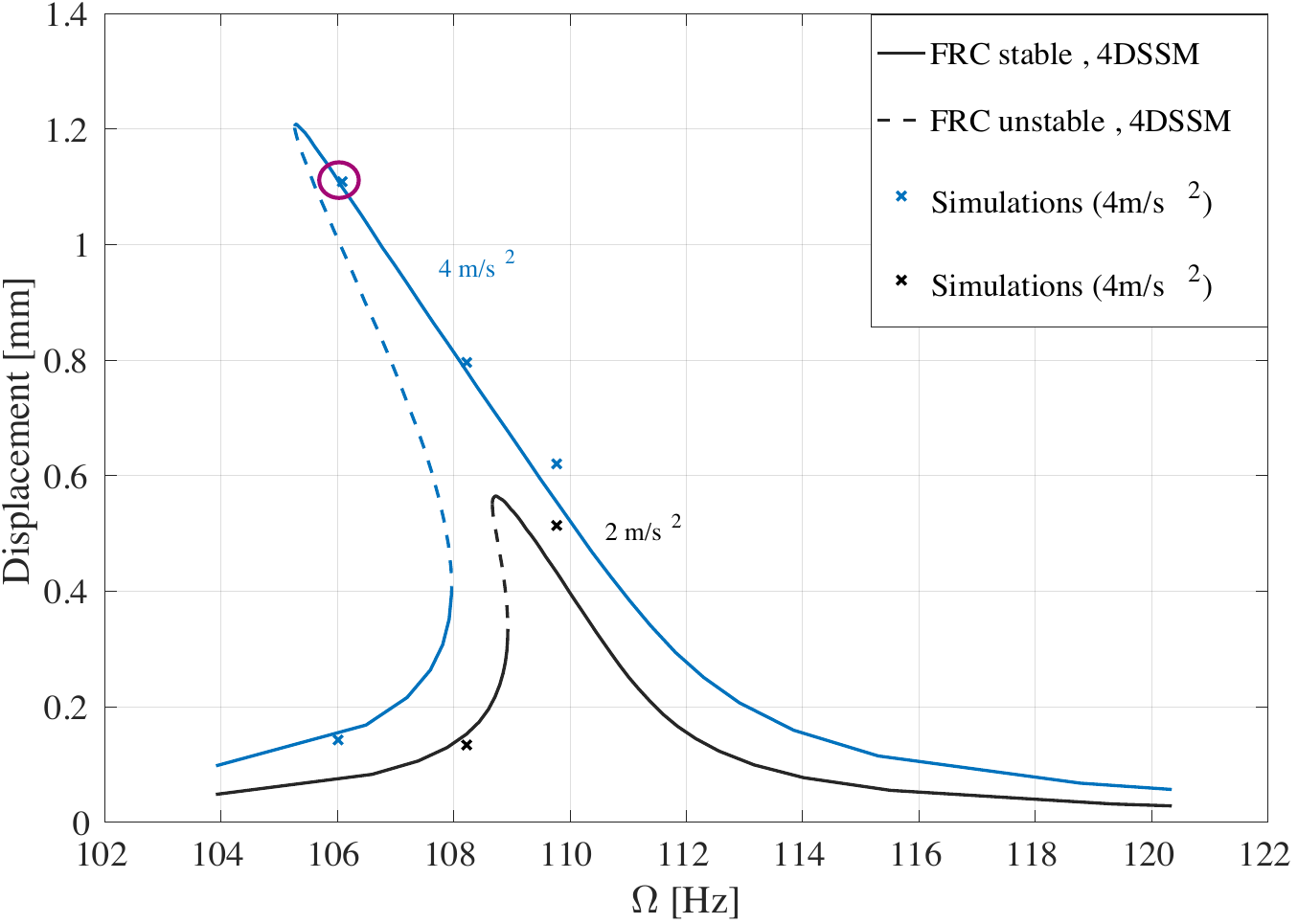}
            \caption{}
            \label{fig:FRFsVsComsol}
        \end{subfigure}\hfill
        \begin{subfigure}{0.5\textwidth}
            \centering
            \includegraphics[width=1\textwidth]{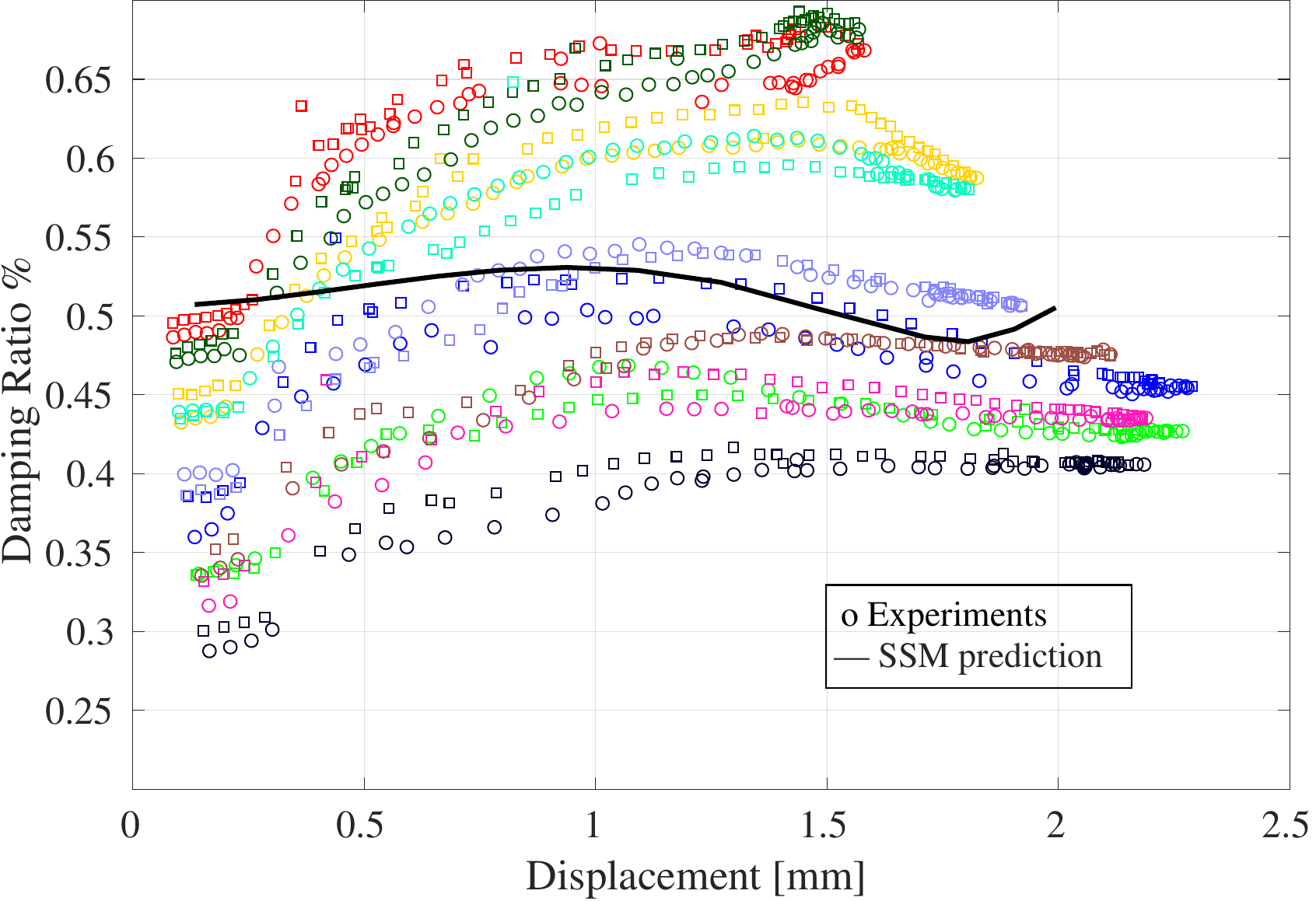}
            \caption{}
            \label{fig:ModalDamping}
        \end{subfigure}
        \caption{Plot (a) shows validations of the maximal forced response displacements at the center of the panel against full simulations (in cross markers). Snapshots of contact and friction stresses at the encircled simulated and predicted responses are plotted in Fig. \ref{fig:stresses}. Plot (b) shows the comparison of modal damping ratios. For clarity, we plot only some of the experiments. As defined experimentally, the displacement measure here is $\sqrt{2} \cross \text{root mean square of the out-of-plane displacements at the center of the panel across one time period}$.}
    \end{figure}

        \begin{figure}[H]
        \centering
        \begin{subfigure}{0.5\textwidth}
            \centering
         \includegraphics[width=1\textwidth]{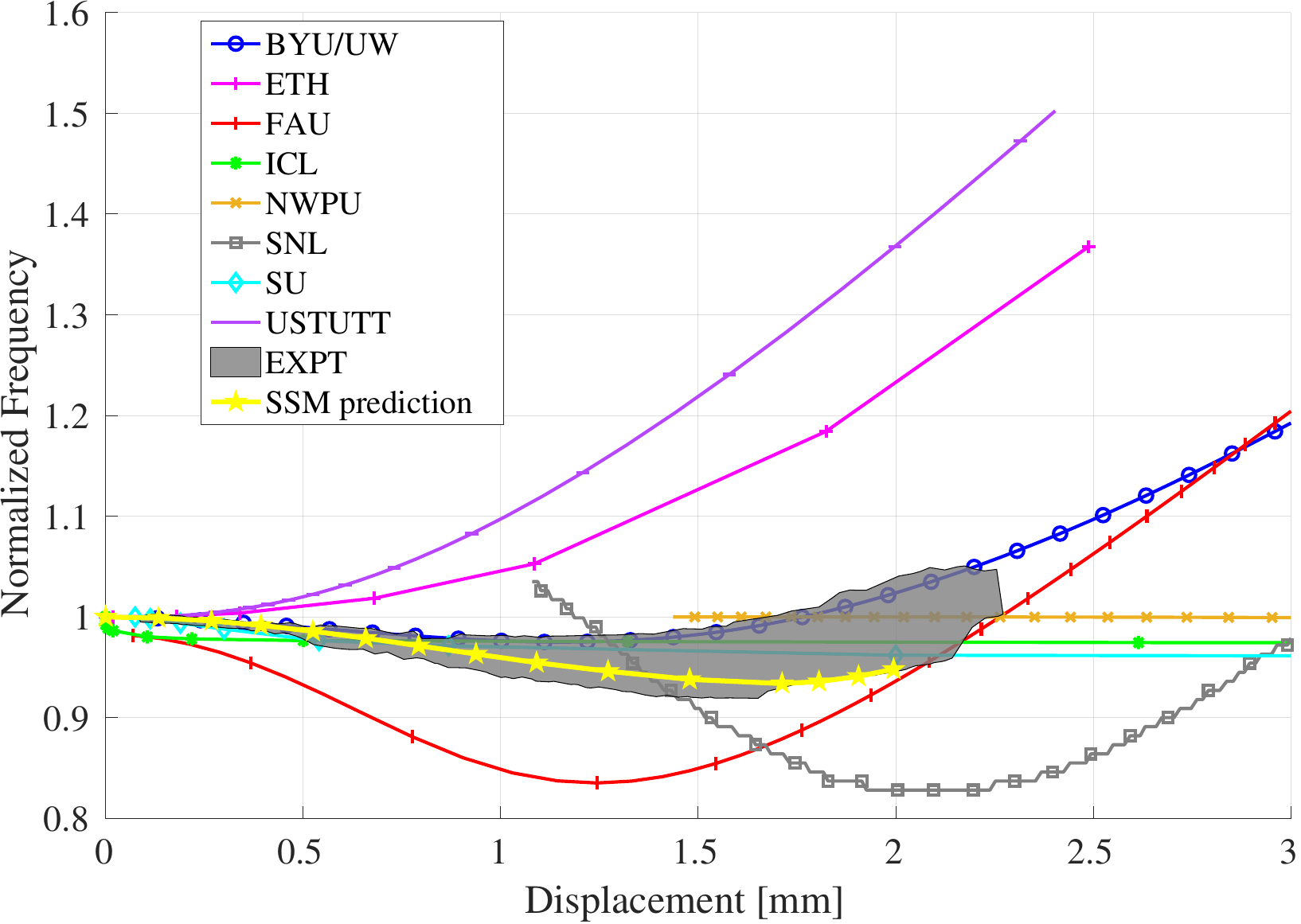}
            \caption{}
      \label{fig:BlindFrequency}
        \end{subfigure}\hfill
        \begin{subfigure}{0.5\textwidth}
            \centering
            \includegraphics[width=1\textwidth]{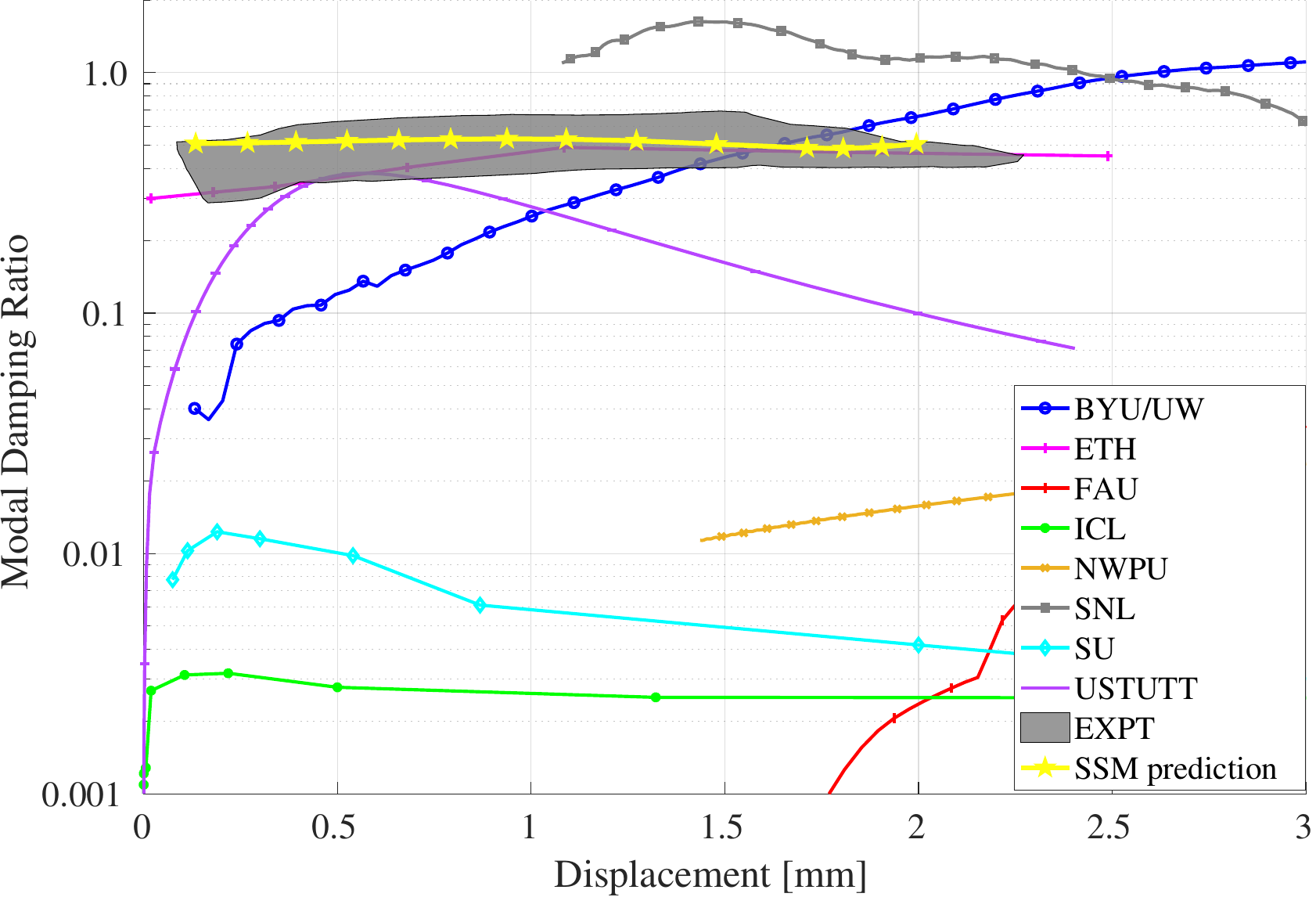}
            \caption{}
            \label{fig:BlindDamping}
        \end{subfigure}
        \caption{Plot (a) compares the predicted backbone curves (in black with star markers) with their experimentally observed range and with blind predicitons results from various approaches presented in \cite{TRCarticle}. Plot (b) compares the corresponding modal damping ratios. For both plots, the displacement measure is $\sqrt{2} \times \text{root mean square of the out-of-plane displacements at the center of the panel across one time period}$.}
    \end{figure}

         \begin{figure}[H]
        \begin{subfigure}{0.5\textwidth}
            \centering
            \includegraphics[width=1\textwidth]{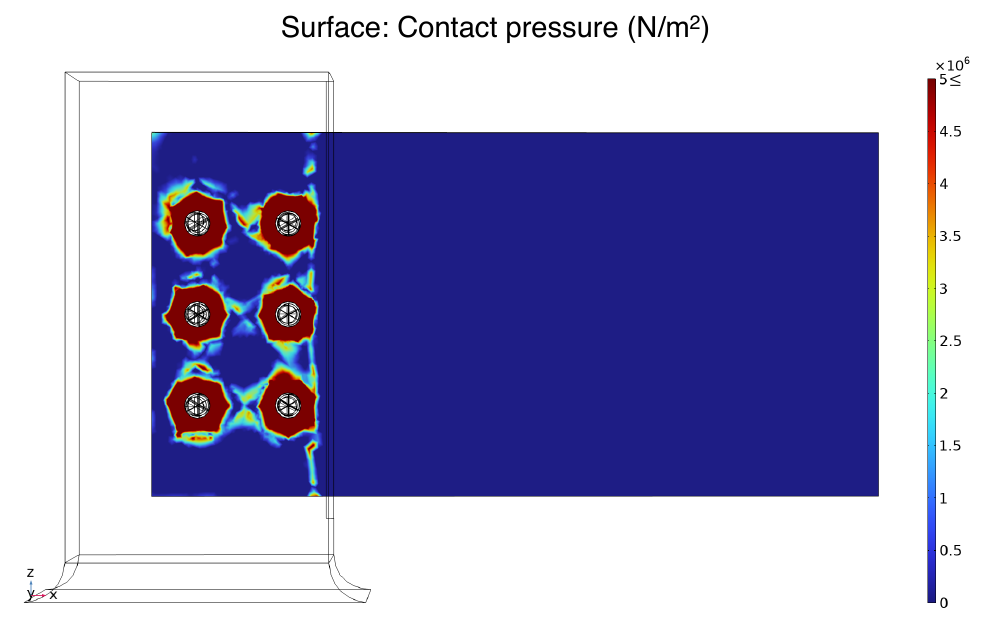}
            \caption{}
            \label{fig:SimPressure}
        \end{subfigure}\hfill
        \begin{subfigure}{0.5\textwidth}
            \centering
            \includegraphics[width=1\textwidth]{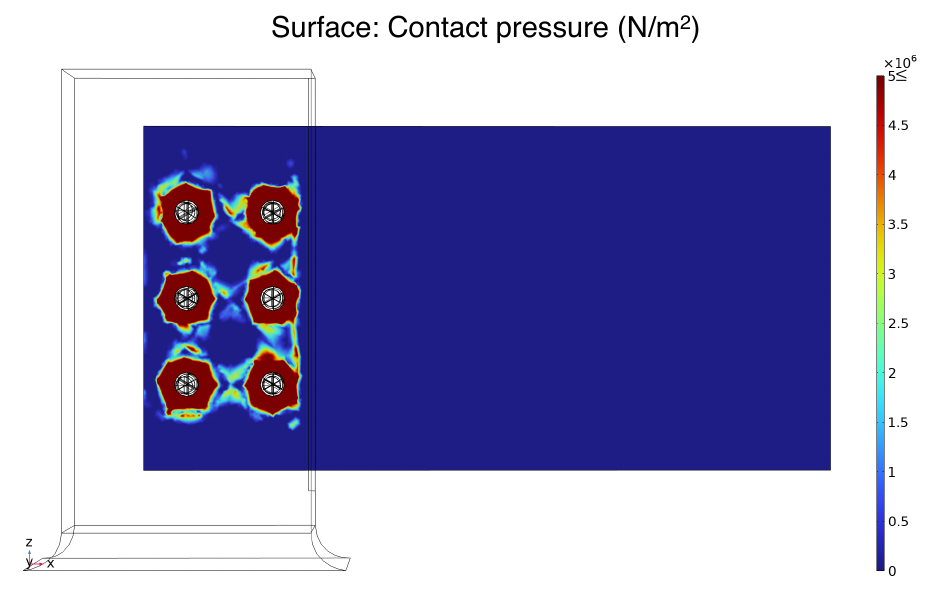}
            \caption{} 
            \label{fig:PredPressure}
        \end{subfigure}
        \vfill         
        \centering
        \begin{subfigure}{0.5\textwidth}
            \centering
         \includegraphics[width=1\textwidth]{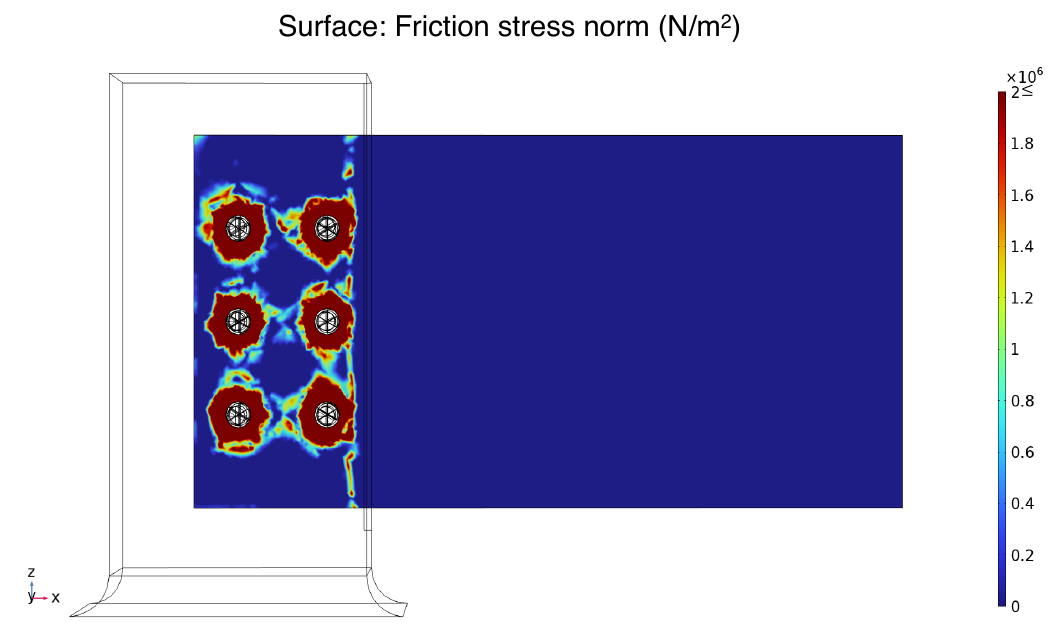}
        \caption{}
            \label{fig:SimFriction}
        \end{subfigure}\hfill
        \begin{subfigure}{0.5\textwidth}
            \centering
            \includegraphics[width=1\textwidth]{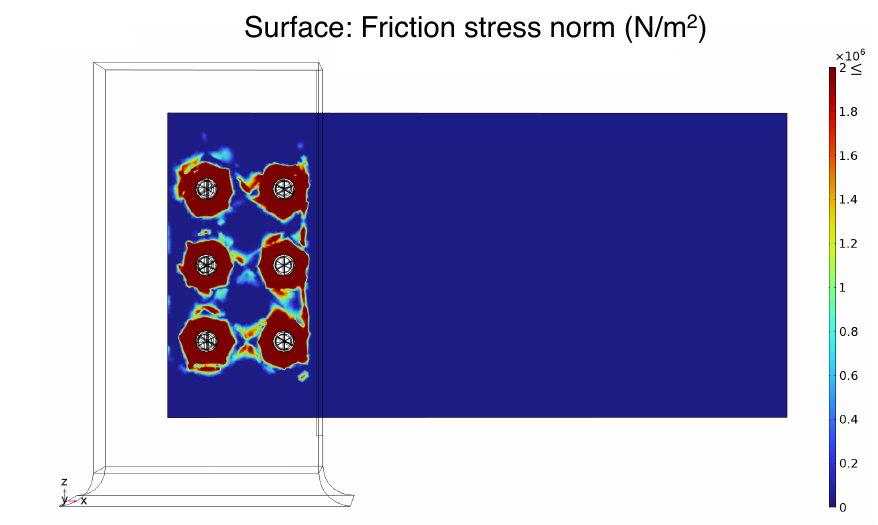}
            \caption{}
          \label{fig:PredFriction}
        \end{subfigure}
        \caption{Plot (a) shows a snapshot of the periodic contact stress distribution on the interface between the panel and the support structure. This is post-processed from the full forced simulation encircled in the sketch in \ref{fig:FRFsVsComsol}. Plot (b) shows the SSM-predicted contact stresses at the same time frame. Plot (c) shows the friction stresses post-processed from the full simulation. Plot (d) shows the SSM-predicted friction stress at the same time frame. For meaningful plots, we upper-bounded the color bars by $2\cross10^{6}$ and $5\cross10^{6}$ $\frac{\text{N}}{\text{m}^{2}}$. Areas close to the bolt holes show higher, localized pressure and friction forces.}
        \label{fig:stresses}
    \end{figure}
        \section{CONCLUSION}
    \label{section:conclusions}

    We have derived a data-driven ROM for a benchmark problem with bolted joints using the recent theory of spectral submanifolds (SSMs). This approach is purely dynamics-based and hence does not require assumptions beyond those needed for accurate, high-fidelity FE modeling of the structure. We use a small number of unforced simulations for training, then construct a ROM by learning the local geometry and the dominant dynamics on a 4D SSM in the phase space of the mechanical system. We then use the ROM to predict the forced responses of the structure and validate them against full simulations and experiments. The SSM-based predictions are accurate not only for the forced response curves, their backbone curves and modal dampings, but also for the local stresses inside the joints. 

    As for limitations, the SSM-based ROM predicts the slow dynamics of the structure which must be activated in the training data. Therefore, to model fast modes, one would need more training trajectories to construct a high-dimensional SSM. The other limitation is the computational bottleneck concerning the simulation times required for the training trajectories. While running full simulations is costly, this dynamics-based approach, however, offers an unrivaled opportunity for accurate prediction of forced responses that are unfeasible using full simulations. Moreover, the forcing type can be arbitrary, as long as it is moderate either in amplitude or speed \cite{Haller2024}. Finally, as discussed in Section \ref{section:TRCBenchmark}, the computational cost of generating a small number of training trajectories can be largely alleviated by using explicit time integration schemes. Another limitation of deterministic ROMs is that joints are often subject to influential parametric uncertainties such as geometry and friction parameters. To explore such a parameter space using our ROM, more training data corresponding to different parameter values will be needed.

    Finally, we have presented predictions of the the forced response of the benchmark structure under periodic forcing. However, our ROM can also be used to make predictions of the response in presence of quasi-periodic forcing, aperiodic forcing \cite{Haller2024}, or even random forcing \cite{xu2024}. Furthermore, the influence of parameter perturbations can be studied by computing the sensitivities of SSM-based ROMs \cite{Li2024}. 
    \section{ACKNOWLEDGEMENTS}
		Paolo Tiso gratefully acknowledges the financial support of the Swiss National Science Foundation for the project “\textbf{Meso}-scale modeling of \textbf{Fri}ction in reduced non-linear interface \textbf{Dy}namics: \textbf{MesoFriDy}".
        \section{Appendix A1}
    \label{section:appendix}
    The autonomous extended normal form of the 4D SSM-reduced dynamics of the TRCBenchmark structure studied in Section \ref{section:TRCBenchmark} reads

\begin{equation}
  \begin{split}
  \dot{\rho}_{1}\rho^{-1}_{1} = &-3.5066-8.788 \rho^{2}_{1}-75.7952 \rho^{2}_{2}\\ &+118.6605 \rho^{4}_{1}-1597.9282 \rho^{2}_{1} \rho^{2}_{2}+1066.6574 \rho^{4}_{2}\\ &-555.7177 \rho^{6}_{1}-28320.9235 \rho^{4}_{1} \rho^{2}_{2}+84941.887 \rho^{2}_{1} \rho^{4}_{2}\\ &-4514.4489 \rho^{6}_{2}\\ &+\mathrm{Re}((-19.7594+27.6132i) \rho_{2} e^{i(-2 \theta_{1}+ \theta_{2})})\\ &+\mathrm{Re}((-659.28713+987.46728i) \rho^{2}_{1} \rho_{2} e^{i(+2 \theta_{1}- \theta_{2})})\\ &+\mathrm{Re}((1415.9495+144.32986i) \rho^{2}_{1} \rho_{2} e^{i(-2 \theta_{1}+ \theta_{2})})\\ &+\mathrm{Re}((59.03333+142.002i) \rho^{3}_{2} e^{i(-2 \theta_{1}+ \theta_{2})})\\ &+\mathrm{Re}((-651.084937-12060.7023i) \rho^{4}_{1} \rho_{2} e^{i(+2 \theta_{1}- \theta_{2})})\\ &+\mathrm{Re}((218.4112+621.5515i) \rho^{2}_{1} \rho^{3}_{2} e^{i(+2 \theta_{1}- \theta_{2})})\\ &+\mathrm{Re}((359.33975-4100.1714i) \rho^{4}_{1} \rho_{2} e^{i(-2 \theta_{1}+ \theta_{2})})\\ &+\mathrm{Re}((172.76657-3258.697i) \rho^{2}_{1} \rho^{3}_{2} e^{i(-2 \theta_{1}+ \theta_{2})})\\ &+\mathrm{Re}((-1870.9639+759.22803i) \rho^{5}_{2} e^{i(-2 \theta_{1}+ \theta_{2})}),\\ 
  \end{split}
  \label{eq:NF_damping1st}
      \end{equation}
  \begin{equation}
  \begin{split}
  \dot{\rho}_{2}\rho^{-1}_{2} = &-5.8897-19.0618 \rho^{2}_{1}-32.9414 \rho^{2}_{2}\\ &+192.3758 \rho^{4}_{1}+586.5103 \rho^{2}_{1} \rho^{2}_{2}+649.996 \rho^{4}_{2}\\ &+118.0441 \rho^{6}_{1}+40945.5356 \rho^{4}_{1} \rho^{2}_{2}-35470.7189 \rho^{2}_{1} \rho^{4}_{2}\\ &+2619.8531 \rho^{6}_{2}\\ &+\mathrm{Re}((-0.00632367+11.6161i) \rho^{2}_{1} \rho^{-1}_{2} e^{i(+2 \theta_{1}- \theta_{2})})\\ &+\mathrm{Re}((0.1532545-156.3585i) \rho^{4}_{1} \rho^{-1}_{2} e^{i(+2 \theta_{1}- \theta_{2})})\\ &+\mathrm{Re}((-394.8948-83.35704i) \rho^{2}_{1} \rho_{2} e^{i(+2 \theta_{1}- \theta_{2})})\\ &+\mathrm{Re}((71.85692+96.46516i) \rho^{2}_{1} \rho_{2} e^{i(-2 \theta_{1}+ \theta_{2})})\\ &+\mathrm{Re}((-135.7648+870.6976i) \rho^{6}_{1} \rho^{-1}_{2} e^{i(+2 \theta_{1}- \theta_{2})})\\ &+\mathrm{Re}((11409.3614+36480.4635i) \rho^{4}_{1} \rho_{2} e^{i(+2 \theta_{1}- \theta_{2})})\\ &+\mathrm{Re}((-4257.7698+6128.2715i) \rho^{2}_{1} \rho^{3}_{2} e^{i(+2 \theta_{1}- \theta_{2})})\\ &+\mathrm{Re}((1957.8524+3153.608i) \rho^{4}_{1} \rho_{2} e^{i(-2 \theta_{1}+ \theta_{2})})\\ &+\mathrm{Re}((3190.849+2672.2477i) \rho^{2}_{1} \rho^{3}_{2} e^{i(-2 \theta_{1}+ \theta_{2})}),
  \end{split}
  \end{equation}
  \begin{equation}
    \begin{split}
   \dot{\theta}_{1} = &+692.9433-904.692 \rho^{2}_{1}-287.2357 \rho^{2}_{2}\\ &+4296.951 \rho^{4}_{1}+12334.1853 \rho^{2}_{1} \rho^{2}_{2}-3410.1834 \rho^{4}_{2}\\ &+1542.461 \rho^{6}_{1}-426040.1353 \rho^{4}_{1} \rho^{2}_{2}-98666.2783 \rho^{2}_{1} \rho^{4}_{2}\\ &+26762.5071 \rho^{6}_{2}\\ &+\mathrm{Im}((-19.7594+27.6132i) \rho_{2} e^{i(-2 \theta_{1}+ \theta_{2})})\\ &+\mathrm{Im}((-659.28713+987.46728i) \rho^{2}_{1} \rho_{2} e^{i(+2 \theta_{1}- \theta_{2})})\\ &+\mathrm{Im}((1415.9495+144.32986i) \rho^{2}_{1} \rho_{2} e^{i(-2 \theta_{1}+ \theta_{2})})\\ &+\mathrm{Im}((59.03333+142.002i) \rho^{3}_{2} e^{i(-2 \theta_{1}+ \theta_{2})})\\ &+\mathrm{Im}((-651.084937-12060.7023i) \rho^{4}_{1} \rho_{2} e^{i(+2 \theta_{1}- \theta_{2})})\\ &+\mathrm{Im}((218.4112+621.5515i) \rho^{2}_{1} \rho^{3}_{2} e^{i(+2 \theta_{1}- \theta_{2})})\\ &+\mathrm{Im}((359.33975-4100.1714i) \rho^{4}_{1} \rho_{2} e^{i(-2 \theta_{1}+ \theta_{2})})\\ &+\mathrm{Im}((172.76657-3258.697i) \rho^{2}_{1} \rho^{3}_{2} e^{i(-2 \theta_{1}+ \theta_{2})})\\ &+\mathrm{Im}((-1870.9639+759.22803i) \rho^{5}_{2} e^{i(-2 \theta_{1}+ \theta_{2})}),\\
   \end{split}
   \end{equation}
   \begin{equation}
   \begin{split}       
   \dot{\theta}_{2} = &+1231.0618+742.5998 \rho^{2}_{1}+602.4179 \rho^{2}_{2}\\ &-5765.8853 \rho^{4}_{1}+1074.2317 \rho^{2}_{1} \rho^{2}_{2}-5042.001 \rho^{4}_{2}\\ &-2414.5128 \rho^{6}_{1}-489483.7568 \rho^{4}_{1} \rho^{2}_{2}-7152.7742 \rho^{2}_{1} \rho^{4}_{2}\\ &-999.5532 \rho^{6}_{2}\\ &+\mathrm{Im}((-0.00632367+11.6161i) \rho^{2}_{1} \rho^{-1}_{2} e^{i(+2 \theta_{1}- \theta_{2})})\\ &+\mathrm{Im}((0.1532545-156.3585i) \rho^{4}_{1} \rho^{-1}_{2} e^{i(+2 \theta_{1}- \theta_{2})})\\ &+\mathrm{Im}((-394.8948-83.35704i) \rho^{2}_{1} \rho_{2} e^{i(+2 \theta_{1}- \theta_{2})})\\ &+\mathrm{Im}((71.85692+96.46516i) \rho^{2}_{1} \rho_{2} e^{i(-2 \theta_{1}+ \theta_{2})})\\ &+\mathrm{Im}((-135.7648+870.6976i) \rho^{6}_{1} \rho^{-1}_{2} e^{i(+2 \theta_{1}- \theta_{2})})\\ &+\mathrm{Im}((11409.3614+36480.4635i) \rho^{4}_{1} \rho_{2} e^{i(+2 \theta_{1}- \theta_{2})})\\ &+\mathrm{Im}((-4257.7698+6128.2715i) \rho^{2}_{1} \rho^{3}_{2} e^{i(+2 \theta_{1}- \theta_{2})})\\ &+\mathrm{Im}((1957.8524+3153.608i) \rho^{4}_{1} \rho_{2} e^{i(-2 \theta_{1}+ \theta_{2})})\\ &+\mathrm{Im}((3190.849+2672.2477i) \rho^{2}_{1} \rho^{3}_{2} e^{i(-2 \theta_{1}+ \theta_{2})}).
    \end{split}
    \end{equation}

        \section{Appendix A2}
    \label{section:appendix_A2}

    Here, we illustrate the microslip dynamics of our FE model of the TRC structure. In Fig.\ref{fig:gaps}, we plot the temporal evolution of gaps on the interface between the bottom surface of the panel and the cantilever. The plots show the activation of contact nonlinearities. Corresponding to the same time instants, Fig. \ref{fig:frictiondissipation} shows the friction dissipation rate on the interface, which is due to slippage of contact pairs. These results are post-processed from the full forced simulation encircled in the sketch in \ref{fig:FRFsVsComsol}.

\begin{figure}[H]
        \begin{subfigure}{0.5\textwidth}
            \centering
            \includegraphics[width=1\textwidth]{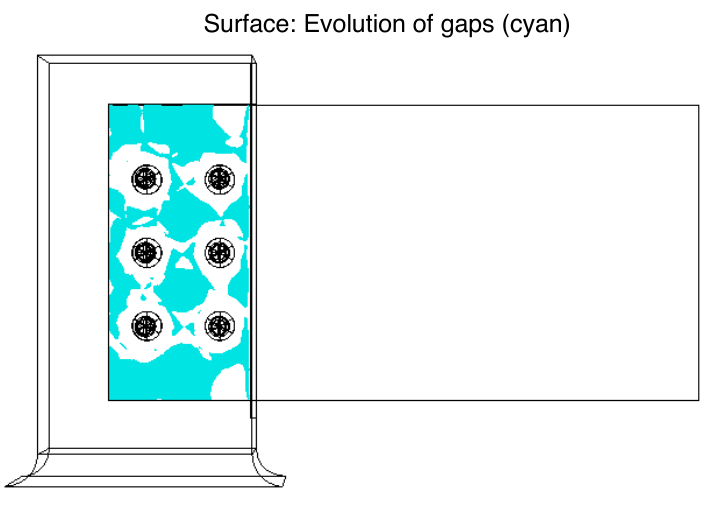}
            \caption{}
            \label{fig:Simcontact90}
        \end{subfigure}\hfill
        \begin{subfigure}{0.5\textwidth}
            \centering
            \includegraphics[width=1\textwidth]{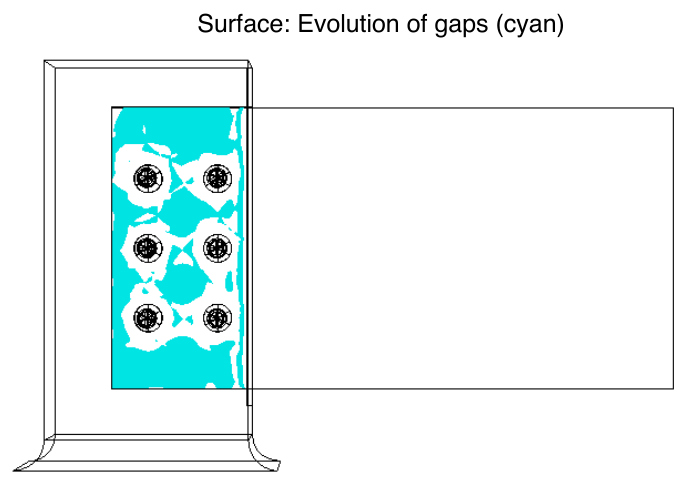}
            \caption{} 
            \label{fig:Simcontact91}
        \end{subfigure}
        \vfill         
        \centering
        \begin{subfigure}{0.5\textwidth}
            \centering
         \includegraphics[width=1\textwidth]{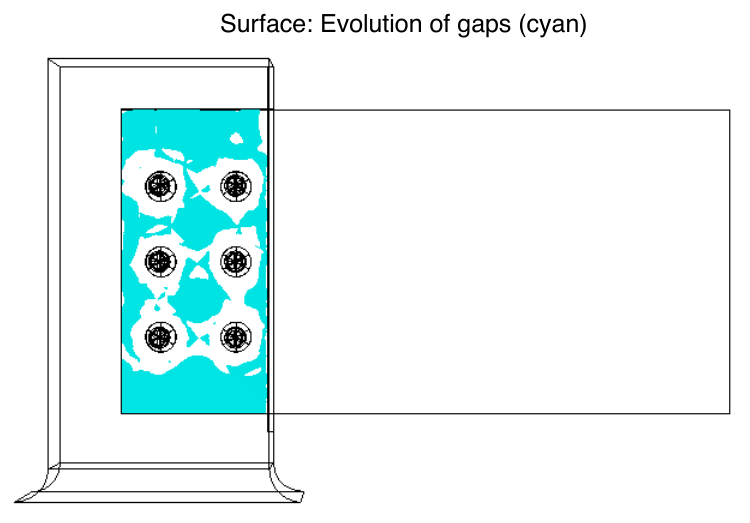}
        \caption{}
            \label{fig:Simcontact92}
        \end{subfigure}\hfill
        \begin{subfigure}{0.5\textwidth}
            \centering
            \includegraphics[width=1\textwidth]{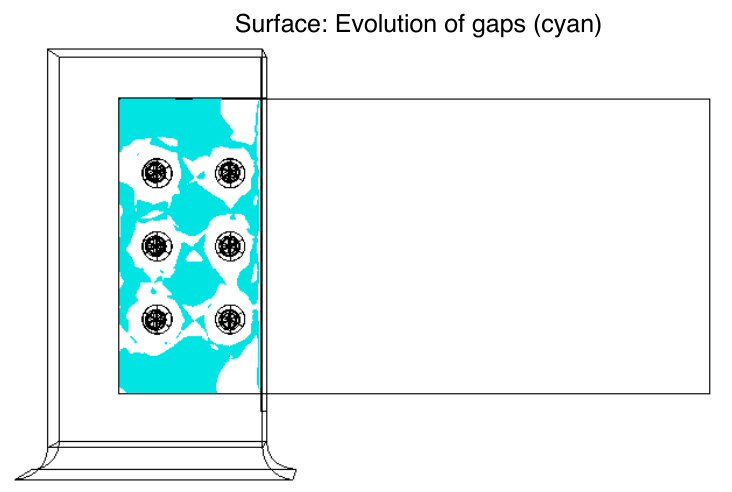}
            \caption{}
          \label{fig:simulationcontact93}
        \end{subfigure}
        \caption{Plots (a,b,c,d) show 4 successive time snapshots illustrating the temporal evolution of the gaps on the bottom surface of the panel. Gaps are indicated in cyan. This is post-processed from the full forced simulation encircled in the sketch in \ref{fig:FRFsVsComsol}.}
        \label{fig:gaps}
    \end{figure} 

    \begin{figure}[H]
        \begin{subfigure}{0.5\textwidth}
            \centering
            \includegraphics[width=1\textwidth]{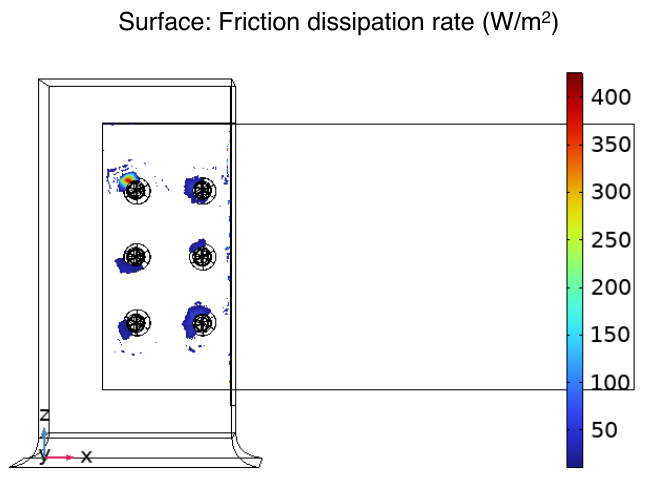}
            \caption{}
            \label{fig:Simfriction90}
        \end{subfigure}\hfill
        \begin{subfigure}{0.5\textwidth}
            \centering
            \includegraphics[width=1\textwidth]{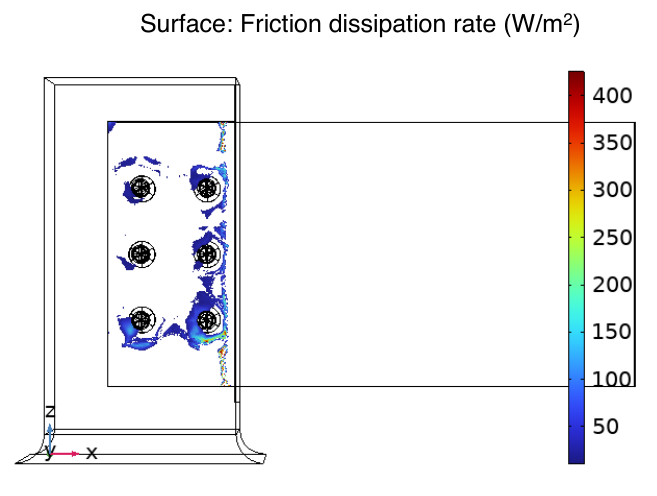}
            \caption{} 
            \label{fig:Simfriction91}
        \end{subfigure}
        \vfill         
        \centering
        \begin{subfigure}{0.5\textwidth}
            \centering
         \includegraphics[width=1\textwidth]{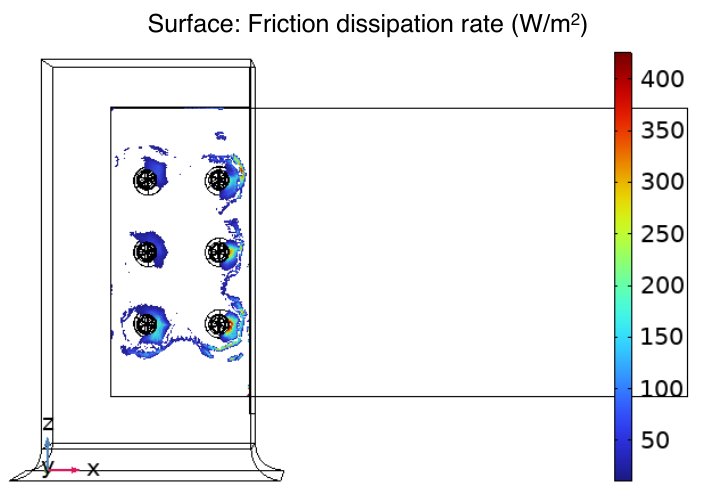}
        \caption{}
            \label{fig:Simfriction92}
        \end{subfigure}\hfill
        \begin{subfigure}{0.5\textwidth}
            \centering
            \includegraphics[width=1\textwidth]{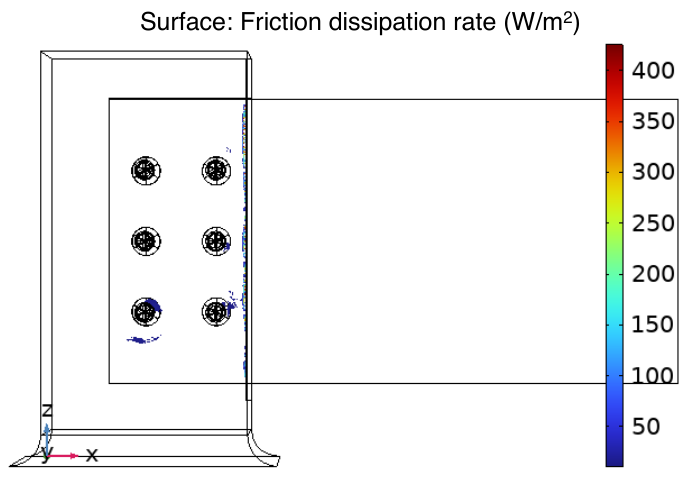}
            \caption{}
          \label{fig:Simfriction93}
        \end{subfigure}
        \caption{Plots (a,b,c,d) show the friction dissipation rate on the bottom surface of the panel at the same 4 successive time snapshots shown in Fig. \ref{fig:gaps}. The friction dissipation rate is computed as the product of the slip force and the slip velocity. This is post-processed from the full forced simulation encircled in the sketch in \ref{fig:FRFsVsComsol}.}
        \label{fig:frictiondissipation}
    \end{figure}

%\section{References}
    \vspace{-1.5ex}
    \bibliographystyle{elsarticle-num}     
    %\renewcommand{\refname}{}
%Where the bibliography will be printed
    \bibliography{references}
  % \printbibliography
\end{document}